\documentclass{amsart}

\usepackage{epsf}
\usepackage[dvips]{graphics}
\usepackage{latexsym}
\usepackage{amssymb}
\usepackage{amsmath}
\usepackage{amscd}
\usepackage{amsthm}

\newtheorem{definition}{Definition}
\newtheorem{propo}{Proposition}
\newtheorem{lemma}{Lemma}
\newtheorem{corollary}{Corollary}
\newtheorem{sublemma}{Sublemma}

\newcommand{\xhomeo}{\mathsf{Homeo}}
\newcommand{\xteich}{\mathsf{Teich}}
\newcommand{\xaut}{\mathsf{Aut}}
\newcommand{\xint}{\mathsf{int}}
\newcommand{\xgraph}{\mathsf{graph}\,}
\newcommand{\xsupp}{\mathsf{supp}\,}
\newcommand{\xiso}{\mathsf{Iso}}
\newcommand{\xtr}{\mathsf{tr}}
\newcommand{\xdegree}{\mathsf{degree}}
\newcommand{\xdev}{\mathsf{dev}\,}
\newcommand{\xexp}{\mathsf{Exp}}
\newcommand{\xid}{\mathsf{id}}
\newcommand{\xpd}{\mathsf{PD}}
\newcommand{\xml}{\mathsf{ML}}
\newcommand{\xhom}{\mathsf{Hom}}
\newcommand{\xcl}{\mathsf{cl}\,}

\newcommand{\xr}{\mathbb{R}}
\newcommand{\xc}{\mathbb{C}}

\newcommand{\xz}{\mathbb{Z}}
\newcommand{\xp}{\mathbb{P}}
\newcommand{\xh}{\mathbb{H}}
\newcommand{\xq}{\mathbb{Q}}

\newcommand{\xrp}{\mathbb{RP}}
\newcommand{\xcp}{\mathbb{CP}}

\newcommand{\xs}{\widetilde{S}}
\newcommand{\xl}{\widetilde{L}}
\newcommand{\xm}{\widetilde{M}}
\newcommand{\xf}{\widetilde{F}}

\newcommand{\iso}{\mathbf{ISO}}
\newcommand{\xo}{\mathbf{O}}
\newcommand{\xpsl}{\mathbf{PSL}}
\newcommand{\xpgl}{\mathbf{PGL}}
\newcommand{\xsl}{\mathbf{SL}}
\newcommand{\xsu}{\mathbf{SU}}
\newcommand{\xso}{\mathbf{SO}}
\newcommand{\xgl}{\mathbf{GL}}
\newcommand{\xO}{\mathbf{O}}
\newcommand{\xoo}{\mathbf{O}_\uparrow}

\newcommand{\xlso}{\mathfrak{so}}
\newcommand{\xliso}{\mathfrak{iso}}

\DeclareMathOperator{\im}{im}

\begin{document}

\title[Lorentz Spacetimes of Constant Curvature]
{Lorentz Spacetimes of Constant Curvature}

\author[Geoffrey Mess]{Geoffrey Mess}

\address{Mathematical Sciences Research Institute\\
Institute for Advanced Study\\
Institut des Hautes \'Etudes Scientifiques\\
and University of California at Los Angeles}

\date{April 19,1990}
\maketitle
\tableofcontents

\section{Introduction}\label{sec1}

A 2+1 dimensional flat spacetime $M$ is a connected 3-manifold with a
flat Lorentz metric, or equivalently a manifold with charts
$f_\alpha:U_\alpha\to\xr^{2+1}$ with transition functions
$f_\alpha\circ f_{\beta}^{-1}$ given by a Lorentz isometry on
$f_{\beta}(U_\alpha\cap U_\beta)$ (assuming $U_\alpha\cap U_\beta$ is
  connected). In the terminology of geometric structures, $M$ has a
  $(\iso(2,1),\xr^{2+1})$-structure. Here $\iso(2,1)$ is the
  automorphism group of $(\xr^{2+1},-dt^2+dx^2+dy^2)$. It is a
  semidirect product $\xo(2,1)\ltimes\xr^{2+1}$ where $\xr^{2+1}$ is
  the subgroup of translations. We restrict attention to oriented
  orthochronous spacetimes; i.e. the transition functions actually lie
  in $\xso(2,1)_0\ltimes\xr^{2+1}$. Any noncompact connected orientable 3-manifold
  immerses in $\xr^{2+1}$ by the immersion theorem \cite{1}, so some
  restriction is needed.

Witten \cite{2} proposes the problem of classifying flat spacetimes
which are topologically of the form $S\times(-1,1)$ and contain
$S\times\{0\}$ as a spacelike hypersurface. Any such spacetime
contains many open sets satisfying the same topological
condition. Thus additional geometrical hypotheses, or else some
natural equivalence relation on spacetimes are necessary for a
reasonable classification. Proposition
~\ref{p6},~\ref{p7},~\ref{p8},~\ref{p15},~\ref{p26} show that it is
natural to consider the class of manifolds which are domains of
dependence (definition~\ref{d4}). Witten also proposes the problems of
classifying spacetimes which are of constant curvature and are
neighbourhoods of a closed spacelike hypersurface, and more generally,
conformally flat Lorentzian spacetimes. The domains of dependence are
spacetimes which are maximal with respect to the property that there is
a closed spacelike hypersurface through each point. In the case of
flat spacetimes, we show in section
~\ref{sec2},~\ref{sec3},~\ref{sec4},~\ref{sec5} that for a given
genus, all domains of dependence form two families, each parametrized
by the cotangent bundle of Teichm\"uller space. All the spacetimes in
one family have an initial singularity and are future complete while
the other family differs only by time reversal.

In the case of anti de~Sitter spacetimes, that is, those which have
constant curvature $-1$, we obtain a complete classification in
section~\ref{sec7}. There is a natural family of maximal domains of
dependence, parametrized by the product of two copies of Teichm\"uller
space, such that given any closed spacelike hypersurface in an anti
de~Sitter spacetime, a neighbourhood of the hypersurface embeds as an open
subset of a unique member of the natural family of spacetimes. These
spacetimes are Lorentzian analogues of the hyperbolic manifolds
corresponding to quasifuchsian groups. There is a convex core and a
consideration of the bending lamination on the two boundary components
of the convex core leads to a new proof of Thurston's theorem on the
existence of a left earthquake and a right earthquake joining any two points in Teichm\"uller space.  In fact this new proof is essentially Thurston's second (and elementary) proof of the existence of earthquakes, interpreted
geometrically in anti de~Sitter space. The spacetimes in this family
have both an initial and a final singularity.  These singularities are qualitatively similar to the initial singularities of the flat spacetimes.

On the other hand, in the case of de~Sitter spacetimes, we only give a
conjectural statement of a complete classification. In
section~\ref{sec6} we construct a family of de~Sitter spacetimes, each
of which is foliated by closed spacelike hypersurfaces, which
conjecturally is the solution to the problem. Each spacetime
corresponds to a projective structure on a Riemann surface at future
infinity. There is another family of spacetimes which have a Riemann
surface at past infinity. Let us remark that most spacetimes in this
family fail to have the property that their universal covers embed in
the model space (that is, de~Sitter space). This is in contrast to the
situation for flat Lorentzian manifolds and for anti de~Sitter manifolds, where the universal covers are convex regions in $2+1$-dimensional Minkowski space or anti de~Sitter space.  Furthermore, in the case of de~Sitter space
there is always an infinite discrete set of spacetimes corresponding to
a given representation of the fundamental group, while in the cases of
curvature $-1$ or 0 the representation determines the spacetime
essentially uniquely.

We have nothing to say about conformally flat Lorentzian manifolds,
nor about the interesting question of foliating the manifolds we
consider by surfaces of constant mean curvature of constant
curvature; however let us mention \cite{52}, cf. section~\ref{sec6}.

Y. Carri\`ere \cite{5} proved a beautiful theorem on the completeness
of closed flat Lorentzian manifolds. In section~\ref{sec8}, we
generalize this to the case of compact Lorentzian manifolds with
spacelike boundary. Thus a flat Lorentzian manifold with boundary is
topologically a product, and is foliated by closed spacelike
hypersurfaces. This result generalizes to Lorentzian manifolds of
constant negative curvature. Probably it also holds in the case of
positive curvature, but we have not shown this. We also make some
remarks about the problem of classifying flat Lorentzian manifolds
containing a closed spacelike hypersurface in the physical case of
3+1 dimensions. While in 2+1 dimensions Einstein's equations require
the spacetime to be flat (or of constant curvature if the cosmological
constant is nonzero) in 3+1 dimensions the spacetimes of constant
curvature are very special, but nonetheless interesting, solutions to
Einstein's equations.

Witten's paper \cite{2} is concerned with the solution of the quantum
theory of general relativity. The definition of the Hilbert space of
the theory depends on knowing the space of classical solutions. In
Witten's theory, classical quantities such as the traces of holonomies
of elements of the fundamental group are observables, but there are
also remarkable new observables which (unlike the holonomy)
distinguish between homotopic, but nonisotopic knots. Other references
to the physical literature are given in \cite{3,4}.

Recall that given a $(G,X)$-manifold $M$, there is a
\textit{development map} $\xdev:\widetilde{M}\to X$ where
we use $\widetilde{\ }$ to denote universal covers. The map $\xdev$ is a submersion and
satisfies $\xdev(\gamma\cdot x)=\rho(\gamma)\xdev(x)$ for each
$\gamma\in\pi_1M$ and a homomorphism $\rho:\pi_1M\to G$
(called the \textit{holonomy map}), determined up to a conjugacy in
$G$. The development $\xdev$ is characterized by the fact that the lift of the
$(G,X)$-structure to $\widetilde{M}$ is given by charts
$\xdev|_{V_\alpha}V_\alpha\to X$.

\textit{Acknowledgments}. I would like to thank Bob Brooks, Jim
Isenberg, Richard Montgomery, Edward Witten, Steve Carlip, Bill
Goldman, Brian Bowditch, Yves Carri\`ere and Fran\c{c}oise Dal'bo for
their interest and encouragement. This work was done at the Mathematical
Sciences Research Institute, the Institut des Hautes \'Etudes Scientifiques, and the Institute for
Advanced Study, and I thank them for their support and
hospitality. This work was partially supported by the NSF.

\section{Fuchsian holonomy}\label{sec2}

Suppose $M$ contains a closed spacelike hypersurface $S$.  At each
point $p\in S$, let $\mathbf{n}$ be the future pointing unit timelike
normal vector. To a unit tangent vector $\mathbf{v}$ at $p$, assign the null
vector $\mathbf{n}+\mathbf{v}$. This identifies the unit tangent
bundle of $S$ with the projectivized bundle of null directions. Let
$q:\iso(2,1)\to\xso(2,1)_0$ be the quotient map. We identify
$\xso(2,1)_0$ with $\xpsl(2,\xr)$.

Then the unit tangent bundle of $S$ is the $\xrp^1$-bundle with structure
group $\xso(2,1)_0$ associated with the flat principal
$\xso(2,1)_0$-bundle over $S$ given by
$$
(\widetilde{S}\times\xso(2,1)_0)/\pi_1S
$$ where $\gamma\in\pi_1S$ acts
on 
$\widetilde{S}\times\xso(2,1)_0$ by $\gamma\cdot(s,g)=(s\cdot\gamma,q(\rho(\gamma))\cdot g)$. So the Euler class of the homomorphism $q\rho:\pi_1S\to\xso(2,1)_0$ is the Euler class of the tangent bundle, that is, $2-2g$ if $S$ is oriented so the positive normal to $S$ is future pointing. If $g>1$, then by Goldman's theorem \cite{6} $q\rho$ maps $\pi_1S$ isomorphically to a discrete cocompact subgroup of $\xso(2,1)_0$. So:

\begin{propo}\label{p1}
If $g>1$ then $q\rho:\pi_1S\to\xso(2,1)_0$ is injective and
$q\rho(\pi_1S)$ is a discrete cocompact subgroup of $\xso(2,1)_0$.
\end{propo}

Proposition~\ref{p1} was conjectured by Witten in \cite{2}. It is possible to avoid Goldman's theorem (cf. the end of section~\ref{sec8}) but we will use an extension of the present argument when considering anti de~Sitter manifolds in section~\ref{sec7}.

For the reader's convenience we give a short exposition of the
Milnor-Wood inequality \cite{7,8} and Goldman's theorem. Suppose
$E\to T$ is an $S^1$-bundle over a closed surface $T$ with structure
group the orientation preserving homeomorphisms of $S^1$ with the
discrete topology. In other words the transition functions for the
bundle are locally constant homeomorphisms of $S^1$. Let $g$ be the
genus of $T$. Regard $T$ as a $4g$-gon $D$, together with $2g$
rectangles (the $1$-handles) attached to represent the generators
$(a_1,a_2,\dots,a_{2g-1},a_{2g})$, and a second $4g$-gon $D'$. Let
$h:\pi_1T\to\xhomeo^+(S^1)$ be the holonomy map. (We regard
$\xhomeo^+(S^1)$ as acting on the \textit{right} of $S^1$.) Regard
$S^1$ as $\xr/\xz$, and consider $H=\xhomeo^+(S^1)$ as the quotient of
the subgroup $H^+\subset\xhomeo^+(\xr)$ defined by $f(x+1)=f(x)+1$ by
its center, the group generated by $z(x)=x+1$. For each generator,
choose a lift $h^\ast(a_i)$ to $H^+$, such that $0\leq0\cdot
h^\ast(a_i)<1$ and $h^\ast(a_{i}^{-1})=(h^\ast(a_i))^{-1}$.

Then there is a section $s$ of the $S^1$-bundle over $D\cup\text{(the
  $1$-handles)}$ which (trivializing the bundle over $D$) equals $0$
and $0\cdot h(a_i)$ over the initial and final sides of the $i$th
$1$-handle and is locally constant along each $1$-handle. The section
is extended radially over $D$, so $s=0$ at the center of $D$ and is
linear along each radius of
$D$. $\prod_{i=1}^{g}[h^\ast(a_{2i-1}),h^\ast(a_{2i})](x)=x-e$ for
  some integer $e$. Choose a constant section say $s'$ over the disc
  $D'$.  The Euler class of the bundle $E$ is the obstruction to extending $s$ over $D'$.  It equals $e\in\xz\cong H^2(T,\xz)$, choosing the generator so
  $\prod_{i=1}^{g}[a_{2i-1},a_{2i}]$ represents the oriented
  boundary of $D\cup\text{(the $1$-handles)}$. Then $0\cdot
  h^\ast(a_1)<1,\; 0\cdot h^\ast(a_1a_2)<2,\; 0\cdot
  h^\ast(a_1a_2a_{1}^{-1})<2,\; 0\cdot
  h^\ast(a_1a_2a_{1}^{-1}a_{2}^{-1})<2,\;  0\cdot h^\ast(a_1a_2a_{1}^{-1}a_{2}^{-1}a_3)<3$ and so on.

It follows that $-e\leq 2g-1$. Similarly $-e\geq1-2g$. These
inequalities can be improved: The genus of a $k$-sheeted covering
space of $T$ is $g'=k(g-1)+1$. The Euler number of the pullback of the
bundle to the covering space is $k\cdot e$. So $|e|\leq 2g-2$, if
$g>1$, and $e=0$ if $g=1$. All values in this range are attained; the
value $2-2g$ is attained by the flat $S^1$-bundle with structure group
$\xpsl(2,\xr)$ (acting on $S^1=\xrp^1$) over a hyperbolic surface
obtained by identifying the unit tangent vectors at each point in the
hyperbolic plane with the circle at infinity. (Of course this
$S^1$-bundle admits no flat structure with structure group $S^1$.)
Benz\'ecri's theorem \cite{9,10} follows: A closed surface of negative
Euler characteristic admits no flat linear connection (not necessarily
torsion-free), and therefore no affine structure. For the Euler number
of the projective tangent bundle is $4-4g$. More generally, an
$\xr^2$-bundle with a flat  $\xgl(2,\xr)$ connection has Euler number
at most $g-1$ in absolute value.

Now consider the components of the space of
representations of a surface group $\pi_1T$, where $T$ has
genus $g>1$, in $\xpsl(2,\xr)$. A representation $\rho$ defines a flat
$S^1$-bundle over $T$ with structure group $\xpsl(2,\xr)$ and Euler
number denoted by $e(\rho)$. If $\rho$ is discrete and faithful this
is the tangent bundle (up to change of orientation) so
$|e|=2g-2$. Conversely Goldman shows in \cite{6} that if
$|e(\rho)|=2g-2$ then $\rho$ is discrete and faithful, so
$e(\rho)=2g-2$ defines the component of the space of representations
which consists of discrete and faithful representations such that the
orientation of $T$ agrees with that of $\xh^2/\rho(\pi_1T)$ and so can
be identified with $\xpsl(2,\xr)\times\xteich(T)$ where $\xteich(T)$ is
Teichm\"uller space. Suppose $\rho$ is a representation with
$e(\rho)=2-2g$. If $\rho$ is not irreducible, it is solvable from
which it would follow that $e(\rho)=0$. So assume that $\rho$ is
irreducible. There is then a $6g-3$-dimensional neighbourhood of $\rho$
in the space of homomorphisms and no function $\xtr^2(\rho(\gamma))$ is
locally constant on the neighbourhood. Suppose $\gamma\in\pi_1T-\{1\}$ and
$\rho(\gamma)$ is an elliptic element of $\xpsl(2,\xr)$, possibly the
identity. By a small change in $\rho$, which doesn't change the Euler
number, we attain that $\rho(\gamma)$ has finite order say $n$. By
\cite{11}, there is a finite sheeted cover $T'$ of $T$ on which $\gamma^n$
is represented by a simple closed curve $C$. We have $e(B')=\chi(T')$
where $B'$ is the flat $S^1$-bundle on $T'$ determined by
$\rho|_{\pi_1(T')}$. But surgery on $C$ gives a 3-manifold $M$ such
that $\partial M=T'\cup T''$ where $\chi(T'')=\chi(T')+2$ and the
bundle extends to $M$ because the holonomy of $C$ is trivial. $T''$,
with the orientation inherited from $M$, is homologous to $-T$. This
contradicts the Milnor-Wood inequality. So $\rho$ is faithful and 
$\im(\rho)$ has no elliptics, and considering the closed subgroups of
$\xpsl(2,\xr)$ (cf. \cite{12}) shows that $\im(\rho)$ is
discrete. R. Brooks has informed me that this proof is given in
\cite{39}. Goldman's proof using foliation is more geometric and the
technique is of much more general use.

A closed spacelike hypersurface $S$ cannot have $g=0$. Otherwise we
would have a spacelike immersion $S\to\xr^{2+1}$. The projection of
$S$ onto $\xr^2$ is then an immersion, which is absurd. We let $g>1$ in the
following proposition, returning to the case $g=1$ below.

\begin{propo}\label{p2}
Given a discrete and faithful $f:\pi_1S\to\xso(2,1)_0$, there is a
flat Lorentz spacetime $M$ with $M=S\times(0,\infty)$ and holonomy
$\rho:\pi_1M\to\iso(2,1)$ equal to $f$.
\end{propo}

\begin{proof}
$f(\pi_1S)$ acts properly discontinuously on each of
$L_+=\{(x,y,t):t^2>x^2+y^2,t>0\}$ and
$L_-=\{(x,y,t):t^2>x^2+y^2,t<0\}$. Take $M$ to be $L_+/f(\pi_1S)$ or
$L_-/f(\pi_1S)$. (This has been well known at least since \cite{13}.)
\end{proof}

Note that $M$ can be chosen either future complete or past
complete. (An orthochronous Lorentz manifold is future complete if every
future pointing timelike or null geodesic segment can be extended to
all positive values of an affine parameter; equivalently, the
developed geodesic extends to a future pointing ray in
$\xr^{2+1}$.  See section~\ref{sec4}.)

\section{Realization of holonomy homomorphisms}\label{sec3}

Fix a homomorphism $f:\pi_1S\to\xso(2,1)_0$. Suppose
$\rho:\pi_1S\to\iso(2,1)$ is a homomorphism with $q\rho=f$. Then
$\rho(\gamma)x=f(\gamma)x+t_\gamma$ for some
$t:G\to\xr^{2+1},\gamma\mapsto t_\gamma$, where $t$ is a 1-cocycle, that is,
$t_{\alpha\beta}=t_\alpha+f(\alpha)t_\beta$. If there exists some $v$
such that $t_\gamma=v-f(\gamma)v$, then $\rho(\gamma)$ differs from $f$ only
by conjugation by the translation $x\mapsto x+v$. The quotient
$H^1(\pi_1S,\xr^{2+1})$ of the cocycles by the coboundaries
corresponds to the space of conjugacy classes of representations
$\rho:\pi_1S\to\iso(2,1)$ such that $q\rho=f$.  We fix an
identification of the quotient space with a subspace of the space of cocycles.

\begin{propo}\label{p3}
Given $f:\pi_1S\to\xso(2,1)_0$ and $M=L_+/f(\pi_1S)$ as in
proposition~\ref{p2}, and any bounded neighbourhood $U$ of $0$ in
$H^1(\pi_1S,\xr^{2+1})$ there exists $C$ such that, letting
$M_C=\{(x,y,t):-t^2+x^2+y^2\leq-C\}/f(\pi_1S)$, there exists a family
of flat Lorentz metrics on $M_C$ parametrized by $U$ such that the
spacetime $M_C(u)$ corresponding to $u\in U$ has holonomy $\rho:\pi_1
M_C\to\iso(2,1)$ such that $q\rho=f$ and
$\rho(\gamma)x=f(\gamma)x+t_\gamma(u)$ where $t_\gamma(u)$ is a
cocycle representing $u\in H^1(\pi_1S,\xr^{2+1})$. (We say that
$M_C(u)$ represents $u$.) Moreover, the spacetimes $M_C(u)$ are future
complete and $S$ is locally strictly convex in each $M_C(u)$.
\end{propo}

\begin{proof}
Consider the spacetime
$M'=\{(x,y,t):t>0,-1\geq-t^2+x^2+y^2\geq-2\}/f(\pi_1S)$. By the
Thurston-Lok holonomy theorem \cite{14,15,16,17}, there is a
neighbourhood $U_0$ of $0$ in $H^1(\pi_1S,\xr^{2+1})$ and a family of
$(\iso(2,1),\xr^{2+1})$-structures on $M'$ parametrized by $U_0$ and
such that the holonomy $\rho:\pi_1(M'(u))\to\iso(2,1)$ represents $u$.

Now choose $K>0$ so that $U\subset K\cdot U_0$. The parametrized
development map $\xdev:\widetilde{M'}\times U_0\to\xr^{2+1}$ satisfies
$$\xdev(\gamma\cdot x,v)=f(\gamma)\cdot x+t_\gamma(v).$$
Now let $M''=K\cdot
M'=\{(x,y,t):t>0,-K\geq-t^2+x^2+y^2\geq-2K\}/f(\pi_1(S))$ and define
$\xdev':\widetilde{M''}\times K\cdot U_0\to\xr^{2+1}$ by
$\xdev'(K\cdot x,K\cdot u)=K\cdot\xdev(x,u)$. Then
$$\xdev'(\gamma x,u)=K\cdot\xdev(\gamma
x/K,u/K)=K(f(\gamma)x/K+t_\gamma(u/K))=f(\gamma)x+t_\gamma(u).$$
So $\xdev'$ is the parametrized developing map of a family of flat
Lorentz structures.  We observe that if $U_0$ is taken sufficiently
small, then $T=\{-t^2+x^2+y^2=-1.5\}/f(\pi_1S)$ remains locally strictly convex and spacelike in the
Lorentz structure determined by any $u\in U_0$. So
$\xdev(\widetilde{T},u)$ is a locally convex complete surface. Later
we will need the following lemma:

\begin{lemma}\label{l1}
Suppose $i:F\to\xr^{2+1}$ is a spacelike immersion of a connected
surface $F$ and $F$ is complete in the induced metric. Then $F$ is the
graph of some function $h:\xr^2\to\xr$, and in particular is embedded.
\end{lemma}

\begin{proof}
The projection $p:F\to\xr^2$ given by
$(i_1(s),i_2(s),i_3(s))\to(i_1(s),i_2(s))$ is a distance increasing
submersion since $i(F)$ is spacelike. Since $F$ is complete, $p$ is a
covering; $F$ is connected so $i(F)$ is a graph.
\end{proof}

By Lemma~\ref{l1}, $S_u=\xdev(\widetilde{T},u)$ is a convex
surface. Given any point $p$ in the future of $S_u$, there is a unique
point $\pi(p)$ on $S_u$ such that the proper time
$\tau(p,q)=\sqrt{-(p-q)\cdot(p-q)}$ is maximized when $q=\pi(p)$. The point $p$
lies on the normal to $S_u$ at $\pi(p)$. $\pi(p)$ is unique and
depends continuously on $p$ because $S_u$ is strictly convex. Now
$p\to(\pi(p),\tau(p,\pi(p)))$ equivariantly identifies the future of
$S_u$ with $S_u\times[0,\infty)$.

So there is a family of flat Lorentz structures on
$\{(x,y,t):-t^2+x^2+y^2<-1\}/f(\pi_1S)$ parametrized by $U$. All are
future complete.
\end{proof}

We have determined all possible holonomies $\rho:\pi_1S\to\iso(2,1)$
and have a family of future complete spacetimes parametrized by this set.

\section{Standard spacetimes}\label{sec4}

\begin{propo}\label{p4}
Let $M$ be a flat Lorentz spacetime. Suppose $M$ contains a closed
spacelike surface $S$ of genus $g>1$. Then $i)$ development embeds
$\widetilde{S}$ in $\xr^{2+1}$ and $ii)$ the coordinate function $t$
is a proper function on $\widetilde{S}$.
\end{propo}

\begin{proof}
$i)$ is by lemma~\ref{l1}. Future pointing unit normals in $\xr^{2+1}$
can be identified with the hyperbolic plane $\xh^2$. Let
$p:\widetilde{S}\to S$ be the projection from the universal covering
space. Let $T:\widetilde{S}\to\xh^2$ be the map such that $T(s)$ is
the future pointing unit normal at $\xdev(s)$. There is an induced map
$T:S\to\xh^2/q\rho(\pi_1S)$. By proposition~\ref{p1},
$\xh^2/q\rho(\pi_1S)$ is a surface of the same genus $g$ as $S$. By
construction the tangent bundle of $S$ is pulled back from 
$\xh^2/q\rho(\pi_1S)$. Since the Euler number of the tangent bundle of
$S$ is $2-2g$, $T:S\to \xh^2/q\rho(\pi_1S)$ has degree $+1$. $T$
induces an isomorphism of fundamental groups so
$T:\widetilde{S}\to\xh^2$ is proper and has degree $1$. Because
$T:\widetilde{S}\to\xh^2$ is proper, the normal $a(s)dt+b(s)dx+c(s)dy$
satisfies $|a(s)|\to\infty$ as $s\to\infty$ on $\widetilde{S}$. Also,
$a^2=b^2+c^2+1$. So
$$dt=-\frac{b(s)}{a(s)}dx-\frac{c(s)}{a(s)}dy.$$
and so $\|dt\|\to1$ as $s\to\infty$, where $\|\quad\|$ denotes the
length in the induced Riemannian metric on
$\widetilde{S}$. Furthermore, since $T$ has degree $1$, for any point
$x\in\xh^2$ and any sufficiently large simple closed curve $C$ in
$\widetilde{S}$, $T(C)$ has winding number $1$ about $x$. Consider
$\widetilde{S}$ as the graph of $t=t(x,y)$, we have that the map
$w:S^1\to S^1$ given by $w(\theta)=(b(s)^2+c(s)^2)^{-\frac{1}{2}}(-b(s),-c(s))$,
where $s=(R\cos\theta,R\sin\theta,t(R\cos\theta,R\sin\theta))$, has
degree $1$ for $R$ sufficiently large. Replace $t$ by a function $u$
which equals $t$ outside a compact set and has exactly one critical
point, which is nondegenerate. This is possible since
$\xdegree(w)=1$. Furthermore $P$ is a maximum or minimum rather than a
saddle, since $\xdegree(w)=1$. Without loss of generality assume $P$
is a minimum. Let $B=u^{-1}\{\epsilon\}$ for some small $\epsilon$. Then
the gradient flow gives a submersion
$F:B\times[0,\infty)\to\widetilde{S}-u^{-1}([0,\epsilon))$. Since
$\|du\|$ is bounded below, $F$ is a 1-1 covering and $|t|\to\infty$ as
$s\to\infty$ on $\widetilde{S}$.
\end{proof}

\begin{corollary}\label{c1}
If a spacetime $M$ contains a closed spacelike surface $S$ of negative Euler
characteristic then $S$ has an isochronous neighbourhood in $M$.
\end{corollary}

\begin{proof}
Otherwise we would have both $t\to\infty$ and $t\to-\infty$ on $\xdev\widetilde{S}$.
\end{proof}

\begin{definition}\label{d1}
A flat Lorentz spacetime $M$ such that $\pi_1M$ is isomorphic to
$\pi_1S$ where $S$ is a closed surface of negative Euler
characteristic is a \textnormal{standard spacetime} if (possibly after a change in
time orientation) $M$ is the quotient of the future in $\xr^{2+1}$ of
a complete spacelike strictly convex surface $\widetilde{S}$ by a
group of Lorentz isometries acting cocompactly on $\widetilde{S}$.
\end{definition}

Given a flat Lorentz spacetime $M$ containing a closed spacelike
surface $S$, we may assume $t\to\infty$ on $\widetilde{S}$ (otherwise
replace $t$ by $-t$). We say $S$ is ``future directed''.

\begin{propo}\label{p5}
Suppose $M$ is a flat Lorentz manifold containing a closed spacelike
future directed surface $S$ of genus $g>1$. There exists a standard
spacetime $M'$ containing a future directed strictly convex surface
$S'$ of genus $g$ and such that $S$ and $S'$ have the same
holonomy. Furthermore $\widetilde{S'}$ can be chosen to lie in the
future of $\widetilde{S}$.
\end{propo}

\begin{proof}
By propositions ~\ref{p1} and ~\ref{p3} there exists a standard
spacetime $M'$ containing a strictly convex surface $S'$ such that
(fixing some identification between $\pi_1S$ and $\pi_1S'$) the
holonomy homomorphisms from $\pi_1S$ and $\pi_1S'$ are equal. Replace
$S'$ by $S'(K)$, the surface of points which lie on the future pointing
normals to $S'$ at proper time $K$. $S$ is compact, so there exists
$K$ sufficiently large such that for each point $p$ in
$\widetilde{S}$ some point of $\widetilde{S'(K)}$ lies in the future of
$p$. Then $\widetilde{S'(K)}$ lies entirely in the future component of
$\xr^{2+1}-\widetilde{S}$, because no timelike line joins two points
in $\widetilde{S'}$.
\end{proof}

\begin{propo}\label{p6}
If $M$ is a spacetime with a closed spacelike surface $S$ of negative Euler characteristic such that
$\pi_1S=\pi_1M$, then there is a spacetime $M''$ containing a
neighbourhood of $S$ in $M$ and containing a standard spacetime $M'$.
\end{propo}

\begin{proof}
We may assume $\widetilde{S}$ is future directed. Let
$\widetilde{S'}$ be a spacelike convex surface such that for each
point $s$ of $\widetilde{S}$ some point of $\widetilde{S'}$ is in the
future of $s$ and $\widetilde{S'}$ is also invariant under
$\rho(\pi_1S)$. $\widetilde{S'}$ together with its future is our
standard spacetime $M'$. It suffices to show that $\rho(\pi_1S)$ acts
properly discontinuously on the region $R$ between the disjoint
spacelike surfaces $\widetilde{S}$ and $\widetilde{S'}$; we then
adjoin $R/\rho(\pi_1S)$ to $M'$ and thicken the resulting manifold
slightly into the past of $S$. (See, for example, \cite{16} for a
formal discussion of thickening a manifold with a geometric
structure.) Define a map $h:\widetilde{S}\times[0,1]\to R$ by
$(s,u)\to(1-u)s+u\phi(s)$ where $\phi(s)$ is the intersection of the
future pointing normal to $\widetilde{S}$ at $s$ with
$\widetilde{S'}$. $h$ is proper: A compact subset $K$ of $R$ lies in a
region $t\leq C$ for some $C$. So $h^{-1}K\subset\widetilde{S}$
contains the set $\{s\in\widetilde{S}:t\leq C\}$. By
proposition~\ref{p4}, $A$ is compact. Now to see $\rho(\pi_1S)$ acts
properly discontinuously on $R$, let $K\subset R$ be compact. Then for
all but finitely many $g\in\rho(\pi_1S)$, $g\cdot h^{-1} K\cap
h^{-1}K=\emptyset=h^{-1}gK\cap h^{-1}K$ so $gK\cap
K=\emptyset$. $\rho(\pi_1S)$ is torsion free, so $R$ covers $R/\rho(\pi_1S)$.
\end{proof}

Since $\rho(\pi_1S)$ acts properly discontinuously on $R$, and $t$ is
a proper function on $R$ which is bounded below, $R/\rho(\pi_1S)$
contains no closed timelike curves. By a standard argument \cite{18},
$R=S\times[0,1]$ where $S\times\{0\}=S,\; S\times\{1\}=S'$, and each
curve $\{s\}\times[0,1]$ is timelike.

We will now consider spacetimes $M$ which are neighbourhoods of
spacelike tori $T$. The linear holonomy
$q\rho:\xz\oplus\xz\to\xso(2,1)_0$ has abelian image. Therefore

\begin{enumerate}

\item
$q\rho$ has image in $\xso(2)$ up to conjugacy, or \label{e1}

\item
the image of $q\rho$ stabilizes a null vector, or \label{e2}

\item	 
the image is hyperbolic, i.e., conjugate into the subgroup of matrices
$$\begin{pmatrix}
1&0&0\\
0&\cosh\lambda&\sinh\lambda\\
0&\sinh\lambda&\cosh\lambda
\end{pmatrix}$$
(in $x,y,t$ coordinates.)\label{e3}
\end{enumerate}

\begin{propo}\label{p7}
In cases ~\ref{e1}) and ~\ref{e2}) the linear holonomy is trivial,
and a neighbourhood of $T$ embeds in a complete spacetime.
\end{propo}

\begin{proof}
1) The holonomy is abelian so lies either in $\xso(2)$ or in the
   subgroup of pure translations. By lemma~\ref{l1}, projection of
   $\xdev\widetilde{T}$ on the plane $t=0$ is bijective. The holonomy
   group acts on the plane $t=0$ and projection intertwines these
   actions. Since the holonomy acts without fixed points on
   $\xdev\widetilde{T}$, the linear holonomy must be trivial.

2) The linear holonomy is a homomorphism
   $\xz\oplus\xz\to\xr\subset\xso(2,1)_0$. If the generators are
   linearly independent over $\xq$, after an arbitrarily small change
   in the linear holonomy, they become dependent and then a basis for
   $\pi_1T$ can be chosen so that $q\rho((0,1))=\xid$. Then if
   $A=\rho((0,1))$, $A\mathbf{x}=\mathbf{x}+\mathbf{n}$. If the linear
   holonomy of $T$ is nontrivial, $\mathbf{n}$ must be a multiple of
   the null vector (say $(1,0,1)$) fixed by the linear holonomy. Let $P$
   be a point on $C=\xdev\widetilde{T}\cap\{y=0\}$. Then
   $A^nP=P+n(1,0,1)$ is on $C$. But the tangent to $C$ is everywhere
   spacelike, so the mean value theorem is contradicted.
\end{proof}

Now suppose the holonomy is purely translational. Since
$\xdev\widetilde{T}$ is spacelike and projects onto $\xr^2$, the
holonomy must map the two generators of $\xz\oplus\xz$ onto
independent spacelike vectors. Choose coordinates so these lie in
$\xr^2$. Then there is a compact spacetime with boundary which
contains a neighbourhood of $T$ in $M$ and is the quotient of 
$C_1\leq
y\leq C_2$ by the holonomy group. Furthermore this spacetime can be
extended to a complete spacetime.

\begin{definition}\label{d2}
A spacetime containing a closed spacelike torus with linear holonomy
nontrivial and lying in the group
$$\{T(\lambda):=\begin{pmatrix}
1&0&0\\
0&\cosh\lambda&\sinh\lambda\\
0&\sinh\lambda&\cosh\lambda
\end{pmatrix},\lambda\in\xr\}$$
is called a \textnormal{standard spacetime} if it is of the form
$\{(x,y,t):t^2>x^2,t>0\}/\langle A,B\rangle$ where
\begin{equation*}
\tag{*}
\begin{aligned}
A\cdot(x,y,t)&=T(\lambda)\cdot(x,y,t)+(0,e,0)\\
B\cdot(x,y,t)&=T(\mu)\cdot(x,y,t)+(0,f,0)
\end{aligned}
\end{equation*}
and $(\lambda,e),(\mu,f)$ are independent in $\xr^2$.

\end{definition}

\begin{propo}\label{p8}
In case ~\ref{e3}) a neighbourhood of $T$ in $M$ embeds in a standard
spacetime (possibly after a change of time orientation).
\end{propo}

\begin{proof}
The holonomy is of the form ~(*) because $\langle A,B\rangle$ is
abelian. The function $t^2-x^2$ on $\widetilde{T}$ is well defined on
$T$ and therefore bounded on $T$. So we have
$\xdev\widetilde{T}\subset\{-C\leq t^2-x^2\leq
C\}\times(-\infty<y<\infty)$. Also, $\xdev\widetilde{T}$ projects
  onto the $(x,y)$-plane, but the projection of the region $t^2-x^2<0$
  is not connected. So $\xdev\widetilde{T}$ must contain some points
  with $t^2-x^2\geq0$. Since $\xdev\widetilde{T}$ is spacelike, it
  does not lie entirely in $t^2-x^2=0$. Suppose without loss of
  generality that $P\in\xdev\widetilde{T}$ and $t(P)>0$,
  $b^2=t^2(P)-x^2(P)>0$. Then we will show $t>0$, $t^2(P)-x^2(P)>0$
  everywhere on $\xdev\widetilde{T}$. Suppose $Q\in\xdev\widetilde{T}$
  and $(t^2-x^2)(Q)=-c^2<0$. Then the orbit of $Q$ lies in a half
  space say $x\leq-c$. However $\xdev\widetilde{T}$ projects 1-1 onto
  the surface $U=t^2(P)-x^2(P)=b^2$ and the projection is equivariant
  with respect to the action ~(*) of $\pi_1T$ on
  $U$. So $Q$ cannot exist. If there existed $R\in\xdev\widetilde{T}$
  with $(t^2-x^2)(R)=0$, then after an arbitrarily small isotopy,
  $\xdev\widetilde{T}$ would be a spacelike surface on which
  $(t^2-x^2)$ changed sign. So (up to a change of sign of $t$)
  $\xdev\widetilde{T}$ lies in a region $t>0,\; a^2\leq t^2-x^2\leq
  b^2,\;-\infty<y<\infty$. The holonomy group acts properly
  discontinuously on that region. It follows that a neighbourhood of
  $T$ embeds in a standard spacetime.
\end{proof}

\section{Domains of dependence and geodesic laminations}\label{sec5}

\begin{definition}\label{d3}

Suppose $M$ is a spacetime containing a closed spacelike future
directed surface $S$ of genus $g>1$. The \textnormal{domain of dependence} of
$\xs$ is the region $D=D(\xs)$ in $\xr^{2+1}$
defined by $x\in\xs$ or $x$ is in the future of
$\xs$ or every future directed timelike or null ray through
$x$ meets $\xs$. And similarly if $S$ is past directed.
\end{definition}

Note that the spacetimes of proposition~\ref{p2} and the standard
spacetimes of definition~\ref{d1} are domains of dependence. We remark
that there is no loss in generality in assuming that $\xs$ is future
directed and $M$ contains a closed smooth strictly convex surface
$S'$ in the future of $S$. The domains of dependence of $\widetilde{S'}$
and $\xs$ are equal.

\begin{propo}\label{p9}
$\pi_1S$ acts properly discontinuously on $D(\xs)$.
\end{propo}

\begin{proof}
Without loss of generality $S$ is future directed. Given $x$ in the
past of $\xs$, there is a neighbourhood $U$ of $x$ such that all
future pointing timelike or null rays through $x$ meet $\xs$ in a
compact set. The proper discontinuity of $\pi_1S$ on $D(\xs)$ then
follows from the proper discontinuity of the action on $\xs$.
\end{proof}

\begin{propo}\label{p10}
Suppose $M_1$ is a standard future directed spacetime and $M_2$ is a
standard past directed spacetime with the same holonomy as $M_1$. Then
the closures of the developments of $\widetilde{M_1}$ and 
$\widetilde{M_2}$ are disjoint.
\end{propo}

\begin{proof}
The intersection of the closures of the developments is a convex
set. Because the restriction of $t$ to the intersection is bounded above
and below, the intersection is a
bounded convex set. If it is nonempty, the holonomy must fix its
barycenter. In this case, the holonomy is conjugate into
$\xso(2,1)$. By definition, the development of a standard spacetime
with holonomy in $\xso(2,1)$ must be a proper subset of the region
$L_+$ or $L_-$ of proposition~\ref{p2}.
\end{proof}

The argument extends to the domains of dependence introduced below.

\begin{definition}\label{d4}
Suppose $M$ is a spacetime containing a closed spacelike surface $S$
of genus $g>1$. $M$ is a \textnormal{domain of dependence} if $M$ is the quotient
of $D(\xs)$ by $\pi_1S$.
\end{definition}

Later we will extend the definition in a natural way to spacetimes of
constant curvature. Recall that a null plane is one on which the
Lorentzian structure is degenerate; equivalently it is a plane
containing a null line but no timelike line. A null plane through the
origin is the subspace orthogonal to the null line through the
origin and contained in the given null plane. All null lines in a
given null plane are parallel.

\begin{propo}\label{p11}

a) The boundary $X$ of a domain $D=D(\xs)$ of dependence in
$\xr^{2+1}$ is convex. b) Any timelike linear function on $\xr^{2+1}$
which increases into the future or the past according as $\xs$ is
future or past directed is bounded below on $D$. c) At any boundary point
$p$, all supporting planes are null or spacelike. d) For each $p\in
X$ there is at least one null supporting plane. Moreover, the cone on
the set of normals to supporting planes is convex and the extreme rays
of this cone are null rays. e) Null rays contained in $X$ and
containing $p$ correspond one-to-one to null supporting planes at $p$.
\end{propo}

\begin{proof}
We assume $D=D(\xs)$ where $\xs$ is future directed and convex. $b)$
follows from proposition~\ref{p10}. Suppose $q$ is not in the domain
of dependence $D$. Then some future directed null ray $l$ through $q$ does not meet
$\xs$, because the union of the set of future directed timelike or
null rays through $p$ which meet $\xs$ is convex set. Let $P$ be the
unique null plane through $q$ containing $l$. Among all planes
$P+(0,0,t)$ parallel to $P$ and disjoint from $D$, the maximum $t_0$
is attained. So $q$ lies in the union of the planes
$P+(0,0,t),t\in(-\infty,t_0]$. So the domain of dependence is an
  intersection of open half spaces with boundaries which are null
  planes. This shows that $D$ is convex and every point of $X$ has a
  null supporting plane. In general, if a convex region in a vector
  space $V$ is defined as the intersection of a set of half spaces
  $\{g(x)\geq C_g:g\in G\}$, the dual vectors to the supporting planes
  at a boundary point $p$ form a compact convex set in (the
  projectivization of) $V^\ast$ and the extreme points are elements of
  $G$ which are minimized at $p$. Given a null ray $l$ in $X$ through
  $p\in X$, there is a supporting plane $Q$ parallel to the null plane
  say $P$ through $l$. Since $D$ is future complete, $Q$ cannot be a
  translate of $P$ towards the past. Since $p\in X$, we must have
  $P=Q$. Now suppose $p\in X$ and $P$ is a null supporting plane of
  $D$ at $p$. Let $l$ be the future directed null ray in $P$ through
  $p$. By construction of the null supporting plane, $l$ does not meet
  $\xs$, so $l$ is not in $D$. But $D$ is future complete, so $l$ is
  in the closure of $D$, so $l$ lies in $X$.
\end{proof}

\begin{propo}\label{p12}
Given a closed oriented hyperbolic surface $S$ and a measured geodesic
lamination on $S$, there is a corresponding flat Lorentz manifold $M$
and an embedding of $S$ in $M$ such that $M$ is a domain of dependence
of $S$ and the linear holonomy of $S$ is the map
$H:\pi_1S\to\xso(2,1)_0$ which is the holonomy of $S$ considered as a
hyperbolic surface.
\end{propo}

\begin{proof}
\cite{19} and \cite{20} are recommended as introductions to geodesic
laminations. First we consider the case where the measured geodesic
lamination $L$ is a finite union of simple closed curves. We start
with a manifold $M_0=L_+/\pi_1S$ as in proposition~\ref{p2}. $S$ lies
in $M_0$ as the quotient of the hyperboloid $\widetilde{S}$ of unit timelike
vectors. Geodesics on $\xs$ are the intersection of planes through the
origin with $\xs$. $M_0$ contains a finite set of totally geodesic
surfaces, each of which is the quotient of the positive light cone in
a 2-dimensional Lorentzian spacetime by an infinite cyclic group
which has two null eigendirections. Now construct a manifold with
boundary, say $M_1$ by splitting $M_0$ along these totally geodesic
surfaces. Suppose the curve $C_i\subset L$ has transverse measure
$a_i$. Let $W_{i,+}$ and $W_{i,-}$ be the corresponding boundary
components of $M_1$. We give $W_{i,+}\times[0,a_i]$ the product
Lorentzian structure, so that the second factor is spacelike and
orthogonal to the first factor and has length $a_i$. Then let
$M=M_1\cup\bigcup_i W_{i,+}\times[0,a_i]$ where $W_{i,+}\times\{a_i\}$
is identified with $W_{i,-}$. Then $M$ is a flat Lorentz manifold
with a closed spacelike convex surface (made from pieces of $S$ and
product annuli). There is a natural homotopy equivalence between $M$
and $M_0$. The linear holonomy is unchanged. The element $t$ of
$H^1(\pi_1S,\xr^{2+1})$ which $M$ represents may be described as
follows: Orient $C_i$. Then $t=\sum_{i}a_i\cdot t_i$ where $t_i$ is the
extension by zero of the cohomology class (with support in an open
annular neighbourhood $N(C_i)$) $\xpd(C_i)\otimes\mathbf{e}_i$ where
$\xpd$ is Poincar\'e duality and $\mathbf{e}_i$ is the unit spacelike
vector fixed by the holonomy of $C_i$. Each $\mathbf{e}_i$ is well defined
in the system of local coefficients on $S$. We fix the orientation of
$\mathbf{e}_i$ so that $\mathbf{u},\mathbf{v},\mathbf{e}_i$ is a
positively oriented triple, where $\mathbf{u},\mathbf{v}$ are the
eigenvectors of the holonomy of $C_i$ with corresponding eigenvalues
less than, respectively greater than, $1$. Now we consider the more
general case of a measured lamination such that the support is not a
finite union of closed leaves. $S-L$ is a union of finitely many
regions $A_i$. Let the preimage of $A_i$ in the universal cover $\xs$
be the disjoint union of regions $A_{ij}$. Fix some region $A_{00}$. We
will translate the cone over each remaining region to obtain the
required spacetime. In case the support of the lamination is a finite
union of closed geodesics, this construction will be the same as the
preceding one, but we will look at it in the universal cover. For each
geodesic $l\in\widetilde{L}$, let $\mathbf{e}_l$ be that unit vector
orthogonal to the plane $P_l$ which meets $S$ in $l$ which points away
from $A_{00}$. Now given a path $h:[0,1]\to\xs$ transverse to $L$, let
$\mathbf{e}:[0,1]\to\xr^{2+1}$ be a continuous function such that
$\mathbf{e}(s)=\mathbf{e}_l$ if $h(s)$ is on the leaf $l$. Given a
region $A_{ij}$, let $h$ be a path with $h(0)\in A_{00},h(1)\in
A_{ij}$. Let $d\mu(s)$ denote the measure on $[0,1]$ determined by the
transverse measure on $\widetilde{L}$ and the map $h$. Let
\begin{equation}\label{eq1}
\mathbf{x}_{ij}=\mathbf{x}_h=\int_{0}^{1}\mathbf{e}(s)d\mu(s).
\end{equation}
Any other path with the same endpoints gives the same value, just as
the ordinary transverse measure is not changed by homotopy with
endpoints fixed.

Consider the union $\xs'$ of the closures of the surfaces
$A_{ij}+\mathbf{x}_{ij}$, which we will show to be
disjoint. $\xs'$ is a translate of $\xs$ by a
vector valued function which is discontinuous only across
$\xl_0$, the sublamination of $\xl$ consisting of
leaves of non-zero mass. Such leaves are the covers of closed isolated
geodesics, and $L$ contains no other closed geodesics. For each
geodesic $l$ in $\xs$ which covers a closed geodesic in $S$, there
correspond two geodesics $l'_{+}$ and $l'_{-}$ which differ by a
translation orthogonal to the planes in which they lie. Attach annuli
as in the case of a lamination supported on finitely many simple
closed curves to obtain a spacelike surface $\xs''$. An alternative
description will be useful to show that the surfaces
$A_{ij}+\mathbf{x}_{ij}$ are disjoint. We define a function
$\mathbf{x}$ on $L_+$. First let $\mathbf{x}$ be $0$ in the open cone over
$A_{00}$. Observe that the lamination $\xl$ determines a
transversely measured lamination of $L_+$, in which the leaves are the
cones over the leaves of $\xl$ (so they are totally geodesic
Lorentzian surfaces). We will also use $L$ to refer to this 2-dimensional
lamination for simplicity.

Choose a basepoint $q_0$ in the cone over $A_{00}$ and, for any $p$ in
$L_+$, choose a path $h$ from $q_0$ to $p$ transverse to
$\xl$. Define $\mathbf{x}(p)=\mathbf{x}_h$ where
$\mathbf{x}_h$ is defined by ~\eqref{eq1}. Then $f(p)=p+\mathbf{x}(p)$ defines a map which is continuous
except across $\xl_0$. By invariance of domain, to show that
$f:L_+-\xl_0\to f(L_+-\xl_0)$ is a homeomorphism it suffices to show
that $f$ is one-to-one. Since $\mathbf{x}$ is constant in complementary
regions and on each leaf, in showing that
$p+\mathbf{x}(p)=q+\mathbf{x}(q)$ implies $p=q$ we may assume that
$p-q$ is null. We have a map $f:L_+\to M$ defined by
$f(p)=p+\mathbf{x}(p)$. Suppose $q\in L_+$ and $p-q$ is null and future
directed.

\begin{lemma}\label{l2}
Suppose $P$ and $Q$ are nondegenerate Lorentzian planes through the
origin in $\xr^{2+1}$. Suppose the line of intersection of $P$ and $Q$
is timelike or null. Then $P$ and $Q$ divide the hyperbolic plane (the
hyperboloid of unit timelike future pointing vectors) into three regions
$E,F,G$ where the closure of $E$ meets $P$ but not $Q$, the closure of
$F$ meets $P$ and $Q$, and the closure of $G$ meets $Q$ but not
$P$. If $\mathbf{m},\mathbf{n}$ are normals to $P,Q$ respectively,
oriented so $\mathbf{m}$ points into $E$ and $\mathbf{n}$ points into
$F$, then the scalar product $\langle\mathbf{m},\mathbf{n}\rangle$ is positive.
\end{lemma}

\begin{proof}
Suppose $P$ is the plane orthogonal to $(0,1,0)$. If a plane $Q$ had a
normal $\mathbf{n}$ orthogonal to $(0,1,0)$, $\mathbf{n}$ would lie in
$P$ and being spacelike, would have a timelike orthogonal complement
in $P$. Then $P$ and $Q$ would have timelike intersection. As any pair
of planes satisfying the conditions of the lemma can be deformed
through pairs of planes satisfying the conditions to any other pair,
and the sign of the scalar product does not change we can assume that
$P$ and $Q$ are almost parallel in which case the lemma is clear. 
\end{proof}

By formula ~\eqref{eq1} applied to the ray from $q$ to $p$,
$\mathbf{x}(p)\cdot(p-q)\geq0$. (Equality holds only if the line
crosses no leaves of the lamination.) Furthermore, since the leaves of
$\xl$ are disjoint planes, lemma~\ref{l2} implies that
$\langle\mathbf{x}(q+t(p-q)),\mathbf{x}(q+t(p-q))\rangle$ increases as
  $t$ increases (and the increase is strict unless the path does not
  meet $\xl$). So $\mathbf{x}(p)$ is spacelike or equal to $0$. Now
  $f(p)-q$ is null or spacelike:
$$\langle f(p)-q,f(p)-q\rangle=\langle
  p-q+\mathbf{x}(p),p-q+\mathbf{x}(p)\rangle
=2\langle p-q,\mathbf{x}(p)\rangle+\langle\mathbf{x}(p),\mathbf{x}(p)\rangle\geq0.$$
Equality holds only if $f(p)=p$. This establishes that if $f(p)=f(q)$
  and $q\in A_{00}$ then $p=q$. But if the basepoint is chosen
  anywhere else in $L_+$, $f$ is changed only by a constant
  translation, so $f$ is injective in $L_+$ and, if $L$ has no
  isolated leaves, $f(L_+)$ is the universal cover of our new
  spacetime. If $L$ has closed leaves, product regions must be added
  to $f(L_+)$. We now wish to show that $\xs''$ is equivariant with
  respect to a homomorphism $\rho:\pi_1S\to\iso(2,1)$ such that the
  linear holonomy $q\rho:\pi_1S\to\xso(2,1)$ is the identity map to
  $\pi_1S$. We will construct a cocycle $t$. Choose a lift
  $\widetilde{p}$ of a basepoint $p$ in $A_{00}$. Given
  $\alpha,\beta\in\pi_1S$, choose paths $h_\alpha,h_\beta:[0,1]\to S$
  with $h_\alpha(0)=h_\alpha(1)=h_\beta(0)=h_\beta(1)=p$ which
  represent $\alpha,\beta$ and are transverse to $L$. Let
  $\widetilde{h_\alpha}$ be the lift of $h_\alpha$ beginning at
  $\widetilde{p}$ and ending at $q$ say. Let $\widetilde{h_\beta}$  be
  the lift of $h_\beta$ beginning at $q$ and ending at $r$ say. Define
  $t_\alpha=\mathbf{x}_{\widetilde{h_\alpha}}$, and similarly for any
  element of $\pi_1S$. Representing $\alpha\beta$ by the path
  $\widetilde{h_\alpha}$ followed by $A\cdot \widetilde{h_\beta}$ we
  have (splitting the integral into two terms)
$$t_{\alpha\beta}=t_\alpha+A\cdot t_\beta.$$
Evidently $\xs'$ is equivariant with respect to the action $\alpha
x=A\cdot x+t_\alpha$. The interior of the union of the cones over the
closures of the surfaces $A_{ij}+\mathbf{x}_h$ (where $h$ is a path
from $A_{00}$ to $A_{ij}$) gives a spacetime $\widetilde{M}$ unless
$L$ contains closed leaves. In this case, we add in regions congruent
to $K_+\times[0,a_i]$ where $K_+$ is a positive light cone in
dimension $1+1$, and $a_i$ is the transverse measure of a closed
leaf. $\pi_1S$ acts and the quotient spacetime $M$ evidently is future
complete and contains a closed convex spacelike surface
$S''=\xs''/\pi_1S$ which will be strictly convex if there are no
closed leaves. Finally, let us show that $\widetilde{M}$ and therefore
also $M$ is a domain of dependence as claimed. Intuitively, $\xs$
stretches and so at least as many timelike and null rays hit $\xs''$
as hit $\xs$.

Suppose that $L$ has no isolated leaves so that $f$ is a
homeomorphism. Suppose $q\in L_+$, $p-q$ is null and future directed,
and $p$ is on $\xs$. Recall that $\xs$ is convex. Let the null and
timelike rays from $q$ meet $\xs$ in the disc $D$ with boundary $C$
corresponding to the null rays. Let $R$ be a spacelike plane in the
future of $q$ and let $D'$ and $f(D)'$ be the radial projections from
$q$ of $D$ and $f(D)'$ on $R$. Let $f(C)'$ be the projection of
$f(C)$. Then from $\langle
f(p)-q,p-q\rangle=\langle\mathbf{x}(p),p-q\rangle\geq0$ with equality
only if $\mathbf{x}(p)=0$, and $f(p)=p$ it follows that $f(C)'$ has
winding number $1$ around any point in the interior of $D'$. So
$f(D)'\supseteq D'$. So $q$ is in the domain of dependence of
$f(\xs)=\xs''$. If $L$ contains isolated leaves, the transverse
measure may be smoothed out over a neighbourhood of each isolated leaf
and the resulting function say $f'$ is a homeomorphism and the
foregoing argument applies. We may use any $A_{ij}$ as a reference
domain in which the function $\mathbf{x}=0$, for some function $f''$
(which differs from $f$ by a constant translation;
$f''(r)=f(r)-\mathbf{x}(s)$ (for some $s\in A_{ij}$)) so for every
point $q$ in $L_+$ which is in the closure of some component of
$L_+-L$, $f(q)=q+\mathbf{x}(q)$ is in the domain of dependence of
$\xs''$. Not every point in $L_+$ lies in the closure of a
complementary component. Let $q$ lie on a leaf which is not a boundary
leaf. Then $f$ is a homeomorphism in a neighbourhood of $q$ so $f(q)$
is in the closure of the domain of dependence of $\xs''$. Since $f$
is open, $f(Q)$ must lie in the domain of dependence of $\xs''$. So
$\widetilde{M}$ is contained in the domain of dependence of $\xs'$. By
construction, through each point in the frontier of $\widetilde{M}$
there is a null ray which does not meet $\xs$. 
So $M$ is a domain of dependence.
\end{proof}

\begin{propo}\label{p13}
Given a spacetime $M$ which is the domain of dependence of a closed
spacelike surface $S$ of genus $g>1$ the geometry of the frontier of
the universal cover determines a measured geodesic lamination on $S$
in the hyperbolic structure on $S$ determined by its linear holonomy,
and this correspondence is inverse to that of Proposition~\ref{p12}.
\end{propo}

\begin{proof}
Let $X$ denote the frontier of $D(\tilde{S})$.
Given a point $r$ on $X$, let $N(r)$ be the union of the future
pointing null rays in $X$ based at $r$. It is a closed subset of
$X$. Let $F(r)$ be the convex hull of $N(r)$.

Claim. Suppose $r\ne s$. Then $F(r)$ and $F(s)$ are disjoint unless
$r$ and $s$ lie on a single null ray. Proof of claim: Choose a timelike
direction and let $p$ project $X$ orthogonally onto an orthogonal
spacelike plane which we identify with $\xr^2$. Rays project to
rays. If there is a unique null ray $l$ through $r$ say, then
$F(r)=l$. If $t\in F(r)\cap F(s)$ and $t\ne r$ then if $t$ was
timelike from $s$ the future of $t$ on $l$ would be timelike from
$s$. But that would contradict $l\subset X$. So $t\in N(s)$. Since
there is a unique null line in $X$ through $t$, we must have that
$r$ and $s$ lie on a single null ray. So we assume that each of $N(r)$,
$N(s)$ consists of at least two null rays. There is a pair of extreme
rays, say $p(l_1),p(l_2)$ through $p(r)$ such that one component of
$\xr^2-l_1\cup l_2$ contains $p(s)$ and is disjoint from $p(N(r))$. We
may choose the direction of projection so that $p(r)=0$ and
$p(l_1)\cup p(l_2)$ span the plane $x=0$. Let $l'_1,l'_2$ be the
extreme rays through $s_1$ so $p(r)$ lies in the component of
$\xr^2-(p(l'_1)\cup p(l'_2))$ which is disjoint from $N(s)$. Since
$p(l'_1),p(l'_2)$ don't cross $p(l_1),p(l_2)$, the $x$-components of
the tangent vectors to $p(l'_1),p(l'_2)$ must be nonnegative. It
follows that the plane sector spanned by $l_1,l_2$ does not meet the
plane $x=0$. So (for some small $c>0$) the plane $x=c$ separates $F(r)$
from $F(s)$.  This proves the claim.

Let $D=D(\widetilde{S})$.
For each $y\in D$ the past of $y$ in $\xr^{2+1}$ is foliated by
strictly convex hyperboloids of constant separation from $y$, so there
is a unique point $a(y)$ on $X$ where separation from $y$ is
maximized, and $a(y)$ depends continuously on $y$. The tangent plane
to the hyperboloid $\{p:\langle p-y,p-y\rangle=\langle
a(y)-y,a(y)-y\rangle\}$ is a spacelike support plane at $a(y)$, and
$y$ is on the timelike normal through $a(y)$ of this spacelike
plane. So every point in $D$ is in some set $F(r)$. Let $E$ be the
subset of points in $X$ with at least 3 supporting null planes, let
$F$ be the subset of points with exactly 2 supporting null planes, and
$G$ the subset of points through which pass exactly one null
ray. Define a lamination $L^\ast$ of $\xm$ by $L^\ast=D-\cup_{r\in
  E}\xint F(r)$. Then $L^\ast$ is the disjoint union of $F(r)$ ($r\in
F$) and the faces of the cones $F(r)$ ($r\in E$) and $L^\ast$ is
closed, so $L^\ast$ is a lamination. By construction, $L^\ast$ is
preserved by the action of the deck group $\pi_1S$, and so it projects to a
lamination on the surface $T=\{y: y\in D,\langle
y-a(y),y-a(y)\rangle=-1\}$. $T$ is convex because the future of $T$
is the union of the convex regions $U(p)=\{y\in D:\langle
y-p,y-p\rangle>-1\}$ where $p\in X$. Each point $y$ on $T$ has a
unique tangent plane, because the support plane at $a(y)$ is parallel
to the tangent plane at $y$. A convex surface has a continuously
varying tangent plane if and only if every point has a unique support
plane, so we have a well defined and continuous map $t: T\to \xh^2$
which takes each point to the point on the hyperboloid of unit
timelike vectors with a parallel tangent plane. Note that $t$ extends naturally
to $t:D\to \xr^{2+1}$; map the ray from $a(p)$ through $p$ isometrically to the ray
from $0$ through $t(p)$. We will define a transverse measure on
$L^\ast$ so that with respect to the induced transverse measure on
$L=t(L^\ast)$, $t$ is the inverse of the map $f$ defined in
proposition~\ref{p13}. We note that $i)$ $t(L)$ is a union of totally
geodesic subspaces $ii)$ $t$ is onto because every timelike unit
vector is the normal to some support plane of $D$ $iii)$ $t$ is
isometric on each component of $\cup_{r\in E}\xint(F(r))$ and so
$t(D-L)$ is open. It follows that $t(L)$ is closed, so $t(L)$ is a
geodesic lamination. Define a vector valued function
$\mathbf{x}:\xh^2\to\xr^{2+1}$ by $t(p)-p=\mathbf{x}(p)$. Then
$d\mathbf{x}$ is a vector valued measure supported on $L$ and (because $t$ takes leaves of $L^\ast$ to leaves of $L$) orthogonal to $L$. So
$d\mathbf{x}$ can be identified with a transverse measure to the
geodesic lamination $L$, and this measure is positive (orienting $L$,
as in proposition~\ref{p12}, so that the normal always points away
from some fixed region) because the regions $t^{-1}(F(r))$ ($r\in E$)
are disjoint. This establishes proposition~\ref{p13}.
\end{proof}

We now have, for each hyperbolic structure on $S$, a bijection from
the space $\xml(S)$ of measured geodesic laminations on $S$ to
$H^1(\pi_1S,\xr^{2+1})$. It is clear that the cocycle constructed in
proposition~\ref{p12} depends continuously on the measured geodesic
lamination, and that the map
$[\mathbf{x}]:\xml(S)\times\xteich(S)\to\mathbf{R}^{6g-6}$ defined by
any continuous trivialization of the bundle $H^1(\pi_1S,\xr^{2+1})$ is
continuous. It follows that the space $\xml(S)$ of measured geodesic
laminations is homeomorphic to $\xr^{6g-6}$. This is analogous to
proving the same result by showing that bending laminations
parametrize the space of quasifuchsian deformations of a Fuchsian
group holding the hyperbolic structure fixed on one boundary
component of the convex hull as in \cite{14,20}. As $\xso(2,1)$
modules, $\xr^{2+1}$ and the Lie algebra $\xlso(2,1)$ are canonically
isomorphic (and so also canonically dual). Witten's paper~\cite{2}
contains the result that the Einstein equation $R_{ij}=0$ can be
regarded as the Euler-Lagrange equations for the Chern-Simons form
$\langle A,dA+\frac{2}{3}[A,A]\rangle$ where $A$ is an $\iso(2,1)$
connection and $\langle\;,\;\rangle$ is the invariant bilinear form on
the Lie algebra $\xliso(2,1)$ for which both the radical $\xr^{2+1}$
and any Levi factor $\xlso(2,1)$ are isotropic; this form is
unique up to multiplication by a constant factor and defines a dual
pairing between $\xlso(2,1)$ and $\xr^{2+1}$. We have shown that
(future directed) domains of dependence form a family parametrized by
the bundle of cohomology spaces $V=H^1(\pi_1S,\xr^{2+1})\to\xteich(S)$,
where $\pi_1S$ varies in Teichm\"uller space. Using the dual pairing
and Poincar\'e duality, this identifies the bundle $V$ with the
cotangent bundle to Teichm\"uller space. By a theorem of
Goldman~\cite{21}, the pairing determines a symplectic structure on
$T^\ast\xteich(S)$, and this is the standard symplectic structure on a
cotangent bundle (rather than the symplectic structure induced from
the symplectic structure of Teichm\"uller space). Of course we can
also identify the bundle $V$ with the tangent bundle to Teichm\"uller
space, because the inner product on $\xlso(2,1)$ determines a
symplectic structure on Teichm\"uller space and so an identification
between the tangent bundle and the cotangent bundle.

We now wish to give a more detailed description of the geometry of the
boundary $X$. First we recall that a geodesic lamination on a hyperbolic
surface has measure $0$. It follows that the map $a\circ
t^{-1}:\xh^2\to X$ maps almost all tangent planes to the countable
set $C$.

\begin{propo}\label{p14}
Suppose $X$ contains two parallel null rays $l_1$ and $l_2$. Then $X$
contains two null rays $l_0$ and $l_3$ parallel to $l_1$ and $l_2$
with basepoints $r$ and $q$ such that the seqment $rq$ is spacelike
and lies in $X$ and from each point in $rq$ there is a unique null ray
in $X$ parallel to $l_1$. Moreover every null ray in $X$ parallel to
$l_1$ has basepoint in $rq$, and there are parallel null rays $l_4$
and $l_5$ through $r$ and $q$ such that through every point on $rq$
there is a null ray in $X$ parallel to $l_4$ and all null rays in $X$
parallel to $l_4$ have basepoint in $rq$. $\pi_1S$ contains an
infinite cyclic subgroup $\langle A\rangle$ which is the stabilizer of
$rq$ and preserves each null ray parallel to $l_1$ or to
$l_2$. $\langle A\rangle$ is the fundamental group of an isolated leaf
of the geodesic lamination determined by $X$. $r,q$ are in $C$. 
\end{propo}

\begin{proof}
Suppose $l_1$ is a null ray based at $r$. Then $F(r)\cap T$ contains a
unique geodesic $g_1$ asymptotic to $l_1$. $t(g_1)$ is a geodesic in
the lamination $L=t(L^\ast)$.

\begin{lemma}\label{l3}
The support of the transverse measure on $L$ is all of $L$.
\end{lemma}

\begin{proof}
Otherwise there would be a closed leaf $C_1$ or a spiralling leaf $C_2$
with transverse measure $0$. But then $t^{-1} C_i$ would be a single
geodesic. Since the cone on $t^{-1}C_i$ separates two complementary
components $F(r_1)$ and $F(r_2)$ we must have $r_1=r_2$. But then
$t^{-1}C_i$ is not part of $L^\ast$ so $C_i$ is not part of $L$.
\end{proof}

\begin{lemma}\label{l4}
If $L$ is a measured geodesic lamination on a closed surface for which
$L$ is the full support of the transverse measure, if $m_1,m_2$ are
distinct geodesics which are asymptotic to each other, then $m_1$ and
$m_2$ are both boundary leaves of a complementary region $A_i$, and so
each of $m_1,m_2$ has transverse measure $0$.
\end{lemma}

\begin{proof}
The proof of lemma~4.5 of \cite{19} establishes this.
\end{proof}

Consider all the null rays $\{l_i\}_{i\in I}$ parallel to $l_1$ in
$X$. For each such null ray $l_i$, we have $l_i\subset N(r_i)$ for a unique
$r_i\not\in G$. $F(r_i)\cap T$ contains a unique geodesic $m_i$
asymptotic to $l_i$. $t(g_i)$ is independent of $i$, because if say
$m_1\ne m_2$ then for all $i\in I$, $t(l_i)=m_1\text{ or }m_2$ (by
lemma~\ref{l4}.) Since $t\circ f=\xid$, this contradicts the fact that
each $m_i$ has transverse measure $0$.

Let $B$ be the set of basepoints of the rays parallel to $l_1$. Then
$D\cap\cup_{b\in B}F(b)=t^{-1}(C)$ for some isolated leaf $C$, and
so $B$ is a spacelike line segment of length $\mu(C)$ where $\mu$ is
the transverse measure. Let $r,q$ be the end points of $B$ and
$l_4,l_5$ the null rays parallel to $l_1$ through $r,q$. $r,q$ are in
$E$ (i.e., they have at least 3 supporting null planes) because
$\xint(F(r))$ and $\xint(F(s))$ are distinct complementary
regions. $\pi_1C$ is the infinite cyclic subgroup $\langle A\rangle$
of the statement of the proposition.
\end{proof}

We call the parts of $X$ which are ruled by parallel null rays the
{\em ruled regions} and their union the {\em ruled part}.

\begin{propo}\label{p15}
Given a spacetime $M=D(\xs)$, there does not exist a spacetime $M'$
strictly containing $M$ such that the development of $\xm'$ does not
meet the ruled part of $X$. In particular, if the holonomy of $M$
fixes a point or the corresponding geodesic lamination has no isolated
leaves, $M$ is a maximal spacetime.
\end{propo}

\begin{proof}
We use the notation $D,X$ as in propositions ~\ref{p13} and
~\ref{p14}. Suppose $M'$ is a spacetime strictly containing $M$, with
development $d:\xm'\to\xr^{2+1}$ and let $\xm''$ be the component of
$d^{-1}(D\cup X)$ which contains $\xm\subset\xm'$. First we consider
the case that $M$ has holonomy conjugate into $\xso(2,1)_0$. We may
assume $M=L_+/\pi_1S$.  

If $0\in\xm''$, then given $g\ne1$ and a
geodesic $C=[x,g\cdot x]/\{g\cdot x=x\}$ in the surface $x\cdot x=-1$,
the cone on $C$ lies in $\xm''$ so $C$ is null homotopic. But $C$ has
nontrivial holonomy $g$. Now suppose $0\ne x\in\xm''\cap
X$. There is a neighbourhood $U$ of $x$ in $X\cap\xm''$ such that $U$
is homeomorphic to an open disc. Since $\pi_1S$ does not act properly
discontinuously on $\xcl L_+$, there exist $q, g\cdot q\in  U$  and a path
$[q,g\cdot q]$ in $U$. Moreover, choosing $x$ on the surface $T$,
there is a homotopy from $[x,g\cdot x]$ to $[q,g\cdot q]$ which moves
each point along a straight line. But then the loop $[x,g\cdot
  x]/\{x=g\cdot x\}$ is homotopic into $U$ and we have a contradiction
as before.

Now suppose that the holonomy does not fix a point. By a ``parallel
region'' we mean a subset of $X$ which is maximal with respect to being
a union of parallel null rays.  The union of the parallel
regions is a closed subset of $X$, because it is the intersection with
$X$ of the closure in $D\cup X$ of the preimage by $t$ of the union of
the isolated leaves of $L$, and $t$ is continuous.  There is no
open subset of $X$ where the direction of the null rays is nonconstant
and varies continuously. For if there was, there would be two null
rays $l_1$ and $l_2$ and a sector $A$ of the circle $S^1$ of null
directions bounded by the directions of $l_1$ and $l_2$, such that the
null ray with direction in $A$ is unique and there is a continuous map
$b:A\times[0,\infty]\to X$ such that $b(\{a\}\times[0,\infty])$ is a
null ray in the direction $a$. But $\pi_1 S\cdot A=S^1$ so outside a
compact region, $X$ is foliated by null rays. Then the holonomy fixes
the barycenter of the convex hull of the set of bases of the null rays
in $X$.  Given $x\in U$ where $U$ is disjoint from the ruled
regions, we will find points $q\in E$ arbitrarily near $x$. If there
was a neighbourhood of $x$ in $U$ consisting of points with unique
supporting plane, then $x$ would have a neighbourhood foliated by null
rays. So there are points arbitrarily near $x$ in $E$ or $F$. If
$\mathbf{n}$ is a timelike direction and the supporting plane
$P(\mathbf{n})$ with normal $\mathbf{n}$ meets $X$ in a single point
(which is the case unless $P(\mathbf{n})$ meets $X$  in a line
segment which is the base of a parallel region, using proposition
~\ref{p14}) then the map $\mathbf{m}\mapsto P(\mathbf{m})$ is
continuous at $\mathbf{n}$. Given a point $r$ in $F$ and a timelike
normal $\mathbf{m}$ to a supporting plane at $r$, there is a timelike
normal $\mathbf{u}$ arbitrarily close to $\mathbf{m}$ for which
$P(\mathbf{u})$ is a supporting plane at some point of $E$, because the set of timelike directions normal to supporting planes at $E$ corresponds (using $t$)
to an open dense subset of the hyperbolic plane, namely
$\xh^2-\xl$. So taking $q=P(\mathbf{m})\cap X$, we have seen that $E$
is dense in the complement of the ruled regions of $X$.

Suppose the holonomy does not fix a point and $x\in
U\subset\xm''\cap X$ and $U$ is a simply connected open set disjoint
from the ruled regions. Then there exist infinitely many points of $E$
in $U$. Given $r\in E$, $F(r)$ is one of the complementary components
of the lamination $L^\ast$. Since there are only finitely many
$\pi_1S$ orbits of complementary components, there are points
$q$ and $g\cdot q$ in $U$ with $q\in E,g\ne1$. Choose $x\in D$ and paths
$[x,g\cdot x]$ and $[q,g\cdot q]$ and join them by a straight line
homotopy  $H:[0,1]\times[0,1]\to D\cup X$ such that
$H:\{s\}\times[0,1]\to D\cup X$ is the affine parametrization of the
straight line joining $H(s,0)\in[x,g\cdot x]$ to $H(s,1)\in[q,g\cdot
  q]$. This projects to a homotopy which moves the loop $[x,g\cdot
  x]/\{x=g\cdot x\}$ into the simply connected neighbourhood $U$. This
is a contradiction. 
\end{proof}

\begin{figure}[t!]

\begin{center}
\epsfysize=4in
\epsffile{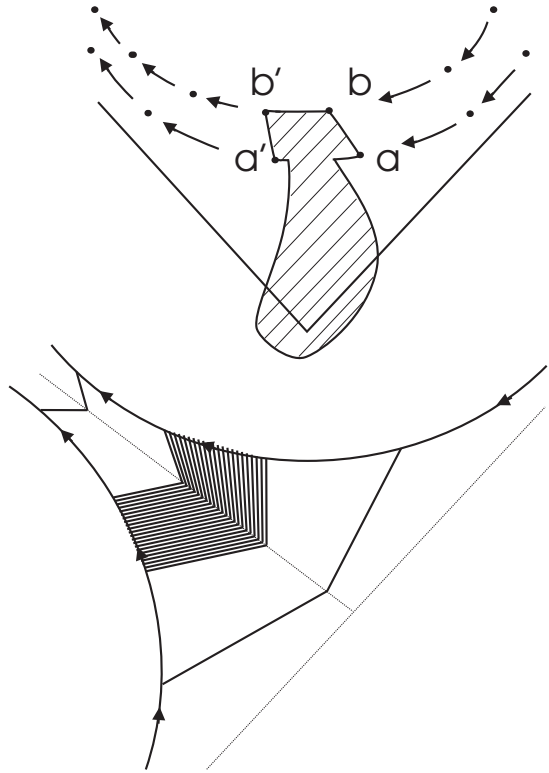}
\end{center}

\caption{\label{fig1}The first example has a fundamental domain which contains the
  fixed point of the holonomy. The development of the universal cover
  is not an embedding. In the second example part of the manifold is
  foliated by closed timelike curves.}
\end{figure}

The need for care is illustrated by Figure~\ref{fig1} which shows the
fundamental domain of flat Lorentz manifolds with infinite cyclic
fundamental group in dimension $1+1$, the development of which contain
the origin or which contain closed timelike curves.

These examples can be ``inserted'' into the ruled part of the
boundary of a domain of dependence, yielding manifolds that contain
closed timelike curves.

In Lorentzian manifolds, spacelike surfaces of constant mean curvature
solve the variational problem of maximizing surface area while
enclosing a given volume. Because spacelike surfaces are locally
graphs with bounded slope, it seems one should be able to deduce the
regularity of maximizers. So one would expect each of the domains of
dependence to have a foliation by surfaces $S_t$ of constant mean curvature such that $t$ is the volume of the past of the surface $S_t$, but we leave both
existence and uniqueness as questions. Moncrief ~\cite{55} discusses
foliations by surfaces of constant mean curvature from a rather
different point of view. One might also ask for a foliation by surfaces
of constant curvature as in ~\cite{52}.

\section{The case of de~Sitter space}\label{sec6}

$\xrp^3$ is separated by the quadric $Q:X^2+Y^2+Z^2=W^2$ into two
regions. The first, $X^2+Y^2+Z^2<W^2$ is the Klein model of hyperbolic
space $\xh^3$. The region $X^2+Y^2+Z^2>W^2$ is the union of the plane
$W=0$ (which can be regarded as the projective plane at infinity
associated with $\xr^3$) and the unbounded component of the complement
of the sphere $x^2+y^2+z^2=1$ in $\xr^3$, identifying $\xr^3$ with the
open set $W\ne0$ in $\xrp^3$ and setting $x=X/W,y=Y/W,z=Z/W$. This
region $X^2+Y^2+Z^2>W^2$ is de~Sitter space. Topologically it is a
twisted $\xr$-bundle over $\xrp^2$. It carries a Lorentzian metric
invariant under the group $\xO(3,1)/\pm1\cong\xoo(3,1)$ of projective
transformations which map the quadric $Q$ to itself. As a homogeneous
space, de~Sitter space is
$\xoo(3,1)/\xoo(2,1)=\xso(3,1)_0/\{\pm1\}\times\xso(2,1)_0$ as the
stabilizer of $(0:0:0:1)$ is $\xoo(2,1)$. Of course
$\xso(3,1)_0\cong\xpsl_2\xc$. The double cover of de~Sitter space is
also called de~Sitter space; it may be realized by the hyperboloid
model, that is, as the set of unit spacelike vectors in $\xr^{3+1}$,
i.e., the vectors $\mathbf{x}$ such that
$\mathbf{x}\cdot\mathbf{x}=1$. De~Sitter space has constant curvature
$+1$. Thus spacelike geodesics focus, and, in the simply connected
de~Sitter space, there is a totally geodesic sphere of constant curvature
$1$ tangent to any spacelike plane in the tangent space of any given
point. In the hyperboloid model these are the intersections of
spacelike planes through the origin with the hyperboloid of unit
spacelike vectors.

We will prefer to work in the Klein model, i.e., in $\xrp^3$. Although
it is not orthochronous, the Klein model of de~Sitter space, like the
Klein model of hyperbolic space, has the advantage that geodesics are
straight lines. A geodesic is timelike, null or spacelike according as
it meets the quadric $Q$ in $2,1,\text{ or }0$ points. Given a point $P$
in de~Sitter space, draw the tangent cone from $P$ to the quadric $Q$;
this is the null cone of $P$. The planes in $\xrp^3$ which do not meet
$Q$ are the totally geodesic spacelike planes.

Recall that a
\textit{projective structure} on a 2-manifold $F$ is a
$(\xpsl_2\xc,\xcp^1)$-structure on $F$, or in other words an atlas of
charts $\phi_\alpha:U_\alpha\to\xcp^1$ where $U_\alpha\subset F$ is
open and connected such that $\phi_\beta\circ\phi_\alpha^{-1}$ is the
restriction of a M\"obius transformation. $F$ has a unique structure
of a complex manifold such that the charts $\phi_\alpha$ are
holomorphic. Classically, the set $P(F)$ of projective structures up
to isotopy on a closed surface $F$ is described as a holomorphic
bundle of affine spaces over the Teichm\"uller space $\xteich(F)$. Any
affine space has an associated vector space of translations, and the
vector bundle over Teichm\"uller space thereby associated with $P(F)$
is the cotangent bundle of Teichm\"uller space. The fiber $\pi^{-1}p$
over $p\in \xteich(F)$ consists of all those projective structures with the
same underlying holomorphic structure (up to isotopy). Given any two
projective structures $A,B\in\pi^{-1}p$, the difference $A-B$ is
defined and is an element of the cotangent space at $p\in
\xteich(F)$. Explicitly, the Schwarzian derivative of one projective
structure with respect to the other defines a holomorphic quadratic
differential on $F$.

Because Teichm\"uller space is a Stein manifold, there must be a
global holomorphic section of $\pi:P(F)\to \xteich(F)$. The Fuchsian family
of projective structures, $f:\xteich(F)\to P(F)$ defined by holomorphically
mapping the universal cover of $F$ to the unit disc, so that the
monodromy is a Fuchsian group, is real but not complex analytic. The
Bers embedding gives a holomorphic section. That is, fix $q\in
\xteich(F)$. For each $p\in \xteich(F)$ there is a quasifuchsian group
$i_p(\pi_1F)\subset\xpsl_2\xc$ with the following properties. The homomorphism $i_p$ is faithful
with discrete image, and the limit set $L_p$ of $i_p(\pi_1F)$ is
homeomorphic to a circle. $L_p$ divides $\xcp^1$ into two regions
$U_+$ and $U_-$. $U_+/i_p(\pi_1F)$ and $U_-/i_p(\pi_1F)$ are Riemann
surfaces corresponding to $p,q$. $U_+/i_p(\pi_1F)$ has a natural
projective structure and this defines a holomorphic section $i:\xteich(F)\to
P(F)$. The holonomy map
$\rho:P(F)\to\xhom(\pi_1F,\xpsl_2\xc)/\xpsl_2\xc$ (where the right
action of $\xpsl_2\xc$ is by conjugation) is a holomorphic submersion.

There is an alternative description of $P(F)$ due to Kulkarni
  ~\cite{34} and Thurston (communicated through ~\cite{35}). This is
  related to work of Apanasov ~\cite{56}. Suppose
  $d:\widetilde{F}\to\xcp^1$ is the development map of a projective
  structure. Over $\xcp^1$ there is an $\xr$-bundle
  $H(\xcp^1)=\xpsl_2\xc/(\xso(2)\ltimes\xc)$ where
  $\xso(2)\ltimes\xc$ is a subgroup of the stabilizer
  $\xc^\ast\ltimes\xc$ of $\infty\in\xcp^1$. The fiber over
  $z\in\xcp^1$ can be thought of as the family of horospheres in
  hyperbolic space tangent to $\xcp^1$ at $z$. There is also the
  bundle $D(\xcp^1)=\xso(3,1)_0/\{\pm1\}\times\xso(2)$ of round discs
  with distinguished point in the interior of the disc: The fiber over
  $z$ is the set of round discs containing $z$ in their interiors. Now
  given $p\in\widetilde{F}$, let $K$ be a subset of $\widetilde{F}$
  containing $p$ such that $d|_K$ is an embedding and the closure of
  $d(K)$ is a round disc and $K$ is maximal with respect to this
  property. Call $K$ a \textit{maximal disc} for $p$. A maximal disc
  $K$ determines a totally geodesic hyperplane in hyperbolic space
  which meets the sphere at infinity in the boundary circle of the
  closure of $d(K)$. Given a maximal disc for $p$, consider larger and
  larger horospheres based at $p$. Eventually there is a horosphere
  tangent to the geodesic hyperplane associated with the maximal
  disc. So a maximal disc has a height, which is an element of the
  fiber of $H(\widetilde{F})$ over $\widetilde{F}$; the height is the
  natural projection $D(\xcp^1)\to H(\xcp^1)$. Of the maximal discs
  for $p$, there is one which is furthest away from $p$, in the sense
  that the height above $p$ is maximal, and this is unique. Indeed, if
  $K_1$ and $K_2$ are maximal discs for $p$, which without loss of 
generality can be
  taken to be the point $\infty\in\xcp^1$, then
  $d(K_1)-\{p\},d(K_2)-\{p\}$ are the exteriors of two circles
  $\partial d(K_1),\partial d(K_2)$ in the Euclidean plane. If
  $\partial d(K_1)$ and $\partial d(K_2)$ are disjoint or meet in a
  single point, then $\widetilde{F}$ contains a region which maps
  conformally to all of $\xcp^1$ or to all but one point of $\xc$,
  which is absurd. If $\partial d(K_1)$ and $\partial d(K_2)$
  intersect in points $A,B$, the outside of the circle with diameter
  $AB$ has preimage $K_3$ and the geodesic hyperplane associated with
  $K_3$ is higher above $p$ than the corresponding hyperplanes for
  $K_1$  or $K_2$. (In the upper half space model with $p$ at
  $\infty$, the hyperplanes are hemispheres or (exceptionally)
  vertical half planes and the higher a hyperplane is above $\infty$
  the lower its tangent horosphere, which is a horizontal plane,
  is). The height defines a section $h:\widetilde{F}\to
  H(\widetilde{F})$ where $H(\widetilde{F})$ is the pullback by
  $d:\widetilde{F}\to\xcp^1$ of the bundle $H(\xcp^1)\to\xcp^1$. The
  highest maximal disc depends continuously on
  $p\in\widetilde{F}$. Define a map
  $D:\widetilde{F}\times[0,\infty)\to\xh^3$ as follows. Given
  $p\in\widetilde{F}$, let $K(p)$ be the highest maximal disc at $p$
  and let $x(p)$ be the point of tangency of $K(p)$ with the horosphere
  based at $p$. Then $D$ maps $\{p\}\times[0,\infty)$ isometrically to the
  geodesic from $x(p)$ to $d(p)\in\xcp^1$, and extends to
  $D:\xf\times[0,\infty]\to\xh^3\cup\xcp^1$ such that
  $D|_{\xf\times\{\infty\}}=d:\xf\to\xcp^1$. $\xf$ may contain two
  overlapping maximal discs $K_1$ and $K_2$, and then $d(K_1\cap K_2)$
  is a crescent shaped part of the sphere, and $K_1\cup K_2$ is
  foliated by curves each of which is mapped by $d$ to a circular arc
  on $\xcp^1$, and these circular arcs all have the same pair of
  endpoints. So $\xf$ contains (at most countably many, but possibly
  $0$) closed regions $B_i$, each topologically a product
  $\xr\times[0,\theta]$ such that each $\xcl(d(\xr\times\{t\}))$ is a
  circular arc joining $p_i,q_i$, and that $d(\xr\times\{t_1\})$ and
  $d(\xr\times\{t_2\})$ intersect at an angle of $|t_1-t_2|$ if
  $|t_1-t_2|<\pi$. Note that $\theta$ may be bigger than $2\pi$. The
  product structure on $B_i$ can be chosen so that each
  $d|_{\{t\}\times[0,\theta]}$ is a submersion with image a circular
  arc orthogonal to each $d(\xr\times\{s\})$.

Now we consider the Lorentzian dual of the Kulkarni-Thurston
construction. Given $p\in\xf$ consider the geodesic through $p$
orthogonal to $K(p)$. Define
$q:(0,\infty)\times\xf\to\xso(3,1)_0/\{\pm1\}\times\xso(2,1)_0$ by
mapping $(t,p)$ to the round disc which bounds a geodesic hyperplane
at distance $t$ from $K(p)$ along the geodesic through $p$ orthogonal
to $K(p)$. Thus this map
$q:(0,\infty)\times\xf\to\xso(3,1)_0/\{\pm1\}\times\xso(2,1)_0$
factors through $\xso(3,1)_0/\{\pm1\}\times\xso(2)$. The map $q$ is locally
injective and so defines a structure of a de~Sitter manifold, i.e. a $(\xso(3,1)_0,\xso(3,1)_0/\{\pm1\}\times\xso(2,1)_0)$-structure on $(0,\infty)\times\xf$. By
construction $q$ is equivariant with respect to the holonomy map.

\begin{propo}\label{p16}
For any genus $g$, there is a family of manifolds, $\pi:X\to
P(F)$, such that each fiber $\pi^{-1}\{a\}$ is a de~Sitter
manifold which is a future complete maximal domain of dependence,
foliated by strictly locally convex spacelike hypersurfaces. Each
admits a conformal compactification at future infinity, and the
projection $\pi$ to $P(F)$ is given by assigning to a given
spacetime the projective structure on the Riemann surface at future
infinity. (There is a second family of past complete spacetimes which
differs only by the covering transformation of the covering of the non
simply connected de~Sitter space by simply connected de~Sitter space.)
\end{propo}

\begin{proof}
For a given projective structure, $q$ defines the structure of a
de~Sitter manifold on $(0,\infty)\times\xf$, and taking the quotient by
$\pi_1F$ represented in $\xpsl_2\xc$ by the holonomy of the projective 
structure gives the required foliated spacetime. For any positive $t$,
the surface $\{t\}\times\xf$ is the dual in de~Sitter space of the
surface in $\xh^3$ at distance $t$ from the locally convex bent
surface $x(\xf)$ in the direction toward infinity. So it is the dual
of a locally strictly convex surface which has $t$ as its smaller
principal curvature at each point. (Given a surface in hyperbolic
space, the dual surface consists of the points dual to the tangent
planes to the surface.) The dual of a strictly convex surface is
strictly convex and spacelike. These spacetimes are maximal domains of
dependence essentially by construction. The rest of the statement of the
proposition is simply descriptive.
\end{proof}

Labourie ~\cite{52} has shown that the spacetimes in proposition
~\ref{p16} are foliated by spacelike convex surfaces of constant
curvature; these surfaces are dual to surfaces of constant curvature
in hyperbolic space. We conjecture that every de~Sitter spacetime
which is a small neighbourhood of a closed oriented spacelike
hypersurface other than a sphere embeds in a unique one of the
spacetimes in proposition ~\ref{p16}. If the hypersurface is a torus
then the projective structure is an affine structure; see~\cite{36}
for a discussion of affine structure on tori. Projective structures on
tori are still classified by the cotangent bundle of Teichm\"uller
space, but the dimension is $2$ rather than $6g-6$ and it is more
natural to think of two affine structures as differing by a
holomorphic differential than by a holomorphic quadratic
differential. It is known that every homomorphism of the fundamental
group of a closed surface of genus $g$ to $\xsl_2\xc$ with image which
is not conjugate into $\xsu(2)$ nor into the stabilizer of a point or
a pair of points in $\xcp^1$ is the holonomy of a projective structure \cite{42}, and that given one projective structure there are infinitely many others with the
same holonomy, but the map from the space of projective structures to
the space of holonomy representations is not a covering map
~\cite{47,48}.

However in the holomorphic description of projective structures, each
fiber $\xc^{3g-3}$ is mapped injectively by the holonomy to
$\xhom(\pi_1F,\xpsl_2\xc)/\xpsl_2\xc$ ~\cite{36}. See also ~\cite{43},
where projective structures whose holonomy is a Fuchsian group are
classified, and  ~\cite{51} where projective structures with holonomy
in $\xpgl(2,\xr)$ are classified. Suppose  we are given a closed
spacelike surface with a neighbourhood modelled on de~Sitter
space. Gallo ~\cite{42} actually shows that given a homomorphism from
the fundamental group of a surface to $\xsl_2\xc$ with image not
conjugate into $\xsu(2)$ nor into the stabilizer of a point or pair of
points in $\xcp^1$, there is a decomposition of the surface into pairs
of pants, such that the holonomy of any pair of pants is
quasifuchsian. After an arbitrarily small change in the holonomy of
each of the boundary circles, the eigenvalue $re^{i\theta}$ of the
holonomy at a fixed point on $\xcp^1$ will have $\theta$ a rational
multiple of $\pi$. Then after passing to a finite cover of the surface
using ~\cite{11}, we may assume that the surface has a decomposition
along curves with purely hyperbolic holonomy into submanifolds $S_i$
whose holonomy groups are quasifuchsian. The bending laminations of
these submanifolds can be approximated by a union of long simple
closed curves, so by another small change in holonomy we can assure
that the submanifolds $S_i$ can be decomposed into submanifolds with
Fuchsian holonomy. So in trying to show that every closed spacelike
surface with a neighbourhood modelled on de~Sitter space lies in a
domain of dependence of a projective structure, one may assume that
the surface has a decomposition into pieces with Fuchsian
holonomy. (There is however a difficulty with this approach which is
pointed out in the proof of proposition ~\ref{p18}.) Even under the
assumption that the holonomy is conjugate into $\xsl(2,\xr)$ it is not
obvious that a closed spacelike surface with a neighbourhood modelled
on de~Sitter space is associated with a projective structure.

Let us
consider the $1+1$ dimensional version of the problem. (Some
information about the action of the automorphism group of de~Sitter
space is given in the proof of the next proposition without being
explicitly stated in the proposition.)

\begin{figure}
\begin{center}
  \rotatebox{-90}{%
      \resizebox{75mm}{!}{\includegraphics{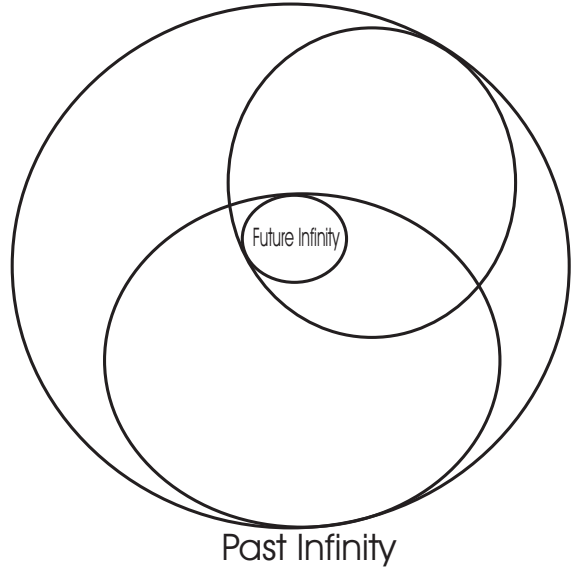}}}
\end{center}
\caption{\label{fig2}$1+1$-dimensional de~Sitter space with some
  typical null geodesics}
\end{figure}

\begin{propo}\label{p17}
Suppose $M$ is a locally de~Sitter $1+1$-dimensional spacetime which
is a small oriented and orthochronous neighbourhood of a closed
spacelike circle $C$. Either

a) $M$ embeds in a complete manifold $M'$ such that the inclusion of
$C$ is a homotopy equivalence and $M'$ has a compactification by a
circle at future infinity and one at past infinity, or else

b) there is a hyperbolic element $\alpha\in\xso(2,1)_0$ such that
$\alpha$ is the holonomy of $C$, and if $P,Q$ are the attracting and
repelling fixed points of $\alpha$ and $R$ the third fixed point of
$\alpha$ ($R$ is the intersection of the tangents at $P$ and $Q$ to the
invariant conic $S^1$ of $\xso(2,1)_0$) there is an arc $PQ$ with
endpoints $P,Q$ on the invariant conic such that the line segments
$RP,RQ$ together with $PQ$ bound an open disc $PQR$ in de~Sitter
space, and $M$ is an open subset of $PQR\cup
PQ/\langle\alpha\rangle$: In short $M$ embeds in a domain of
dependence of a circle with $(\xpsl_2\xr,\xrp^1)$-structure whose
universal cover embeds in $\xrp^1$ and whose holonomy is linear, or

c) $M$ embeds in a domain of dependence $D(M)$ of a circle with
$(\xpsl_2\xr,\xrp^1)$-structure and unipotent holonomy such that the
universal embeds in $\xrp^1$. Then the universal cover
$\widetilde{D(M)}$ of $D(M)$ is the complement in $\xrp^2$ of the
union of the closure $\xcl\xh^2$ and the tangent line to a point 
$P\in S^1=\partial\xcl\xh^2$.
\end{propo}

\begin{proof}

$1+1$-dimensional de~Sitter spacetime $X_{1+1}$ is the complement, in
the double cover $S^2$ of $\xrp^2$, of two antipodal round discs
$D,-D$. Let $\partial X_{1+1}$ be the boundary of $X_{1+1}$: $\partial
X_{1+1}$ is the union of a circle at past infinity and one at future
infinity. Note that $X_{1+1}$ carries two rulings by null geodesics. A
great circle on $S^2$ which is tangent to $D,-D$ at points $P,-P$ in
$S^2$ is divided by $P,-P$ into two open semicircles.  See Figure~\ref{fig2},
in which stereographic projection has been used to draw
$X_{1+1}$ in the plane. Each semicircle is a null geodesic in
de~Sitter space. The null geodesics are divided into two disjoint
families, the left and right rulings of $X_{1+1}$. Recall that if a
manifold is locally modelled on a homogeneous space $G/H$ with
transition functions in $G$, the manifold is also locally modelled on
the universal cover of $G/H$ which is a homogeneous space with
automorphism group in general a proper covering of $G$. The identity
component of the isometry group of the universal cover
$\widetilde{X_{1+1}}$ is the universal cover $\widetilde{\xso(2,1)_0}$
of $\xso(2,1)_0$. The holonomy of $C$ is a well-defined element (up to
conjugacy) $a\in\widetilde{\xso(2,1)_0}$ and acts on the boundary
$\partial\widetilde{X_{1+1}}$. Suppose $a$ acts on one (and therefore
both) components of $\partial\widetilde{X_{1+1}}$ without fixed
points. Then any orbit is unbounded above and below in
$\xr=\partial\widetilde{X_{1+1}}$. Then the projection of
$\xdev\widetilde{C}$ down say the left null ruling to either of the
boundary components must be a submersion onto. Then given any point in
$\widetilde{X_{1+1}}\cup\partial\widetilde{X_{1+1}}$, either the
future pointing or the past pointing timelike and null rays meet
$\xdev\widetilde{C}$ in a compact interval, and we deduce that $a$
acts properly discontinuously on
$\widetilde{X_{1+1}}\cup\partial\widetilde{X_{1+1}}$. Thus $M$ embeds
in $M'=\widetilde{X_{1+1}}/\langle a\rangle$. In particular if the
image of the holonomy is trivial in $\xso(2,1)_0$, then $C$ is
embedded in a finite cover of $X_{1+1}$ and the complete manifold $M'$
is that finite cover.

\begin{figure}

\begin{center}
  \rotatebox{-90}{%
      \resizebox{40mm}{!}{\includegraphics{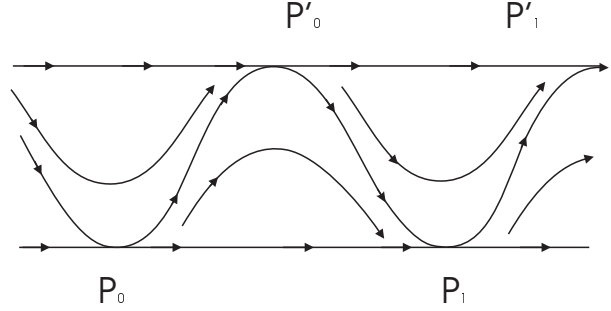}}}
\end{center}
\caption{\label{fig3}The universal cover of $1+1$-dimensional
de~Sitter space and the flow lines of a unipotent subgroup}
\end{figure}

\begin{figure}
\begin{center}
  \rotatebox{-90}{%
      \resizebox{40mm}{!}{\includegraphics{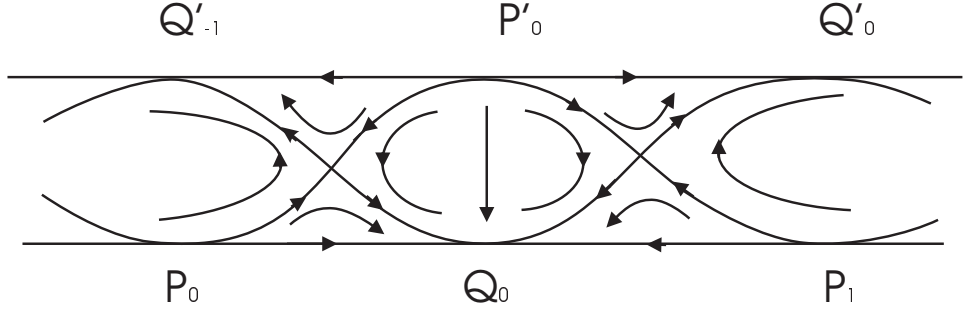}}}
\end{center}

\caption{\label{fig4}The universal cover of $1+1$-dimensional
de~Sitter space and the flow lines of a hyperbolic subgroup}
\end{figure}

On the other hand, if $a$ has a fixed point $P=P_0$ on say the past
component of $\partial\widetilde{X_{1+1}}$, then all the translates
$P_n$, $n\in\xz$ of $P$ by the deck group of the covering
$\partial\widetilde{X_{1+1}}\to \partial X_{1+1}$ are also fixed. The
lifts of the antipodal points are also fixed: These are points $P'_n$
on the future boundary such that there is a null geodesic in the left
ruling joining $P_n$ to $P'_n$ and a null geodesic in the right ruling
joining $P'_n$ to $P_{n+1}$. Now we distinguish two cases. Either $a$
is the lift of a hyperbolic element of $\xso(2,1)_0$ or else $a$ is
the lift of a unipotent element of $\xso(2,1)_0$.  See Figure~\ref{fig3}
and Figure~\ref{fig4}. In the hyperbolic case, there are
fixed points $Q_n,Q'_n$ with $Q_n$ between $P_n$ and $P_{n+1}$, $Q'_n$
between $P'_n$ and $P'_{n+1}$. There are also fixed points $R_n$ on
the intersections of the null geodesics from $Q_n$ to $Q'_n$ and the
null geodesic from $P'_n$ to $P_{n+1}$. Let $p_L:\widetilde{X_{1+1}}\cup \partial\widetilde{X_{1+1}}\to
\partial\widetilde{X_{1+1}}$ be the projection down the left ruling to
the past boundary. Because $\widetilde{C}$ is spacelike
$p_L|_{\widetilde{C}}$ is a homeomorphism onto its image. Since $p_L$
intertwines the action of $\pi_1C$ on $\widetilde{C}$ with the action
of $a$ on the past component of  $\partial\widetilde{X_{1+1}}$,
$p_L(\widetilde{C})$ is an open interval containing no fixed point of
$a$. Therefore $p_L(\widetilde{C})$ is a bounded interval. The
endpoints of the interval $p_L(\widetilde{C})$ are fixed by $a$. So
the interval must be (up to a deck transformation and possibly
interchanging the labels $P,Q$) the interval $P_0,Q_0$. Then the
projection $p_R(\widetilde{C})$ must be one of the intervals
$P_0Q_0,Q_0P_1,P_1Q_1$. In the first of these three cases $C$ lies in
the domain of dependence of the quotient of the interval $P_0Q_0$ by
$\langle a\rangle$; the universal cover of this domain of dependence
is the open triangle $P_0R_0Q_0$, and there is a partial
compactification by adding the quotient of the interval
$P_0Q_0$. Similarly in the third case $C$ lies in a future complete
domain of dependence. The second case cannot arise because any
$\langle a\rangle$ orbit in the open quadrangle $Q_0R_0P'_0R_1$ (whose
boundary consists of null geodesic segments) accumulates at $P'_0$ and
$Q_0$ which is impossible since $\widetilde{C}$ is spacelike.

The unipotent case is similar but a little simpler.
\end{proof}

\begin{propo}\label{p18}
Let $M$ be a de~Sitter manifold which is a small neighbourhood of a
closed spacelike torus. Then (perhaps after a change of time
orientation) there is a torus with affine structure such that $M$
embeds in the spacetime associated to the torus by proposition ~\ref{p16}.
\end{propo}

\begin{proof}
First of all, the holonomy of $M$ cannot be conjugate into the
subgroup $\xsu(2)$ of $\xsl_2\xc$. If it was, then projection onto
$\xcp^1$ along the family of lines which meet at the fixed point
inside $\xh^3$ of the subgroup conjugate to $\xsu(2)$
would define a Riemannian structure on $T$ of constant positive
curvature, which is impossible for a closed surface of positive
genus. Suppose we show that after an arbitrarily small change of the
geometric structure on a neighbourhood of $T$, the resulting manifold
embeds in a domain of dependence of some torus at future infinity
having an affine structure. Then the manifold $M$ is associated to an
affine structure by proposition ~\ref{p16} because the space of all
de~Sitter structures on a germ of a neighbourhood of $T$ is a union of
components each of which is mapped by the holonomy submersively to the
space of conjugacy classes of representations of the fundamental
group, and the map from the space of affine structures to the
universal cover of the space of conjugacy classes of representations
of $\xz\oplus\xz$ which are not conjugate into $\xsu(2)$ is proper, so
if arbitrarily near geometric structures embed in a domain of
dependence then so does the given geometric structure on $M$. However,
the properness may fail in the case of the nonelementary
representations of the fundamental group of a surface of genus greater
than one. Also, if a finite cover of $M$ embeds in a domain of
dependence of an affine structure, $M$ itself does. This applies
equally to the question of surfaces of genus greater than one. We can
then suppose that the holonomy is hyperbolic or loxodromic with
rotational angle a rational multiple of $\pi$, rather than
unipotent. So after passing to a finite cover of the spacetime we may
assume that the holonomy is conjugate into the identity component of
the diagonal subgroup of $\xsl_2\xr$. Then by another small change in
the geometric structure we have (in an appropriate basis) one of the
summands of $\xz\oplus\xz$ is in the kernel of the holonomy. Now an
infinite cyclic cover $T'$ of $T$ lies in the simply connected
de~Sitter space, which is the complement in the three-sphere of two
antipodal round discs. There is an isometric $S^1$ action on de~Sitter
space which commutes with the holonomy of $T$ and fixes an axis, that
is, the great circle which passes through the two fixed points of the
action of the holonomy of $T$ on either of the spheres at infinity of
de~Sitter space. There is a global section to this action so on $T$
there is an $S^1$-valued function say $\theta$ defined everywhere
except at finitely many points which are the images in $T$ of the
intersection of $T'$ with the axis. Euler characteristic
considerations imply that this finite set is empty. 
Now consider the submanifold $T_0$
of $T$ which has a fixed value of $\theta$. Its preimage in $T'$ is
the preimage by the developing map of a plane (i.e., an equatorial
two-sphere) in de~Sitter space through the axis. Each component of
$T_0$ develops into a spacelike curve. Proposition ~\ref{p17} applies
and we conclude that (in the $\xrp^3$ model) $T'$ lies between the
tangent planes to the sphere at infinity at the two points fixed by
the holonomy of $T$, from which it follows readily that $M$ embeds in
a domain of dependence.
\end{proof}

\section{Anti de~Sitter manifolds}\label{sec7}

As with de~Sitter space, we will work with a synthetic description of
anti de~Sitter space. Consider the quadric $Q:AD-BC=0$ in $\xrp^3$. It
divides $\xrp^3$ into two regions, $AD-BC>0$ and $AD-BC<0$. The
quadric is the Segre embedding of $\xrp^1\times\xrp^1$ in $\xrp^3$:
$((X_1:Y_1),(X_2:Y_2))\mapsto(X_1X_2:X_1Y_2:Y_1X_2:Y_1Y_2)$. (The
Segre embedding of $\xp^n\times\xp^m$ in $\xp^{nm+n+m}$ comes from
projectivizing the canonical map $V\times W\to V\otimes W,
(\mathbf{v},\mathbf{w})\mapsto \mathbf{v}\otimes\mathbf{w}$, where
$\xp^n,\xp^m$ are the projectivizations of vector spaces $V,W$.) Thus
there are two families, the left and right rulings, of straight lines
on $Q$. Through each point passes one line of each ruling. Any line of
the left ruling intersects each line of the right ruling in exactly
one point, and vice versa. So the lines of the right ruling are
parametrized by the points of any one line of the left ruling, and
vice versa. Now recall that in $\xr^4=\xr^2\otimes\xr^2$, the
decomposable vectors are those that lie on the quadric $AD=BC$. The
subgroup $G_1$ of $\xgl(\xr^4)$ which preserves the decomposable
vectors is the group which preserves the quadratic form $AD-BC$ up to
scale. $G_1$ has an index two subgroup
$G_2=\xgl(\xr^2)\times\xgl(\xr^2)/\langle(aI,I)=(I,aI)\rangle$;
the non-identity coset of $G_2$ in $G_1$ is generated by an involution
which exchanges the left and right factors of the tensor product
$\xr^4=\xr^2\otimes\xr^2$. $G_2$ contains as a normal subgroup, with
quotient group $\xz/2\oplus\xz/2\oplus\xr_{>0}^{\ast}$ (where
$\xr^{\ast}_{>0}$ is the multiplicative group of the positive reals),
the group $\xso(2,2)_0=\xsl_2\xr\times\xsl_2\xr/(-I,I)=(I,-I)$. So the
identity component of the automorphism group of the quadric $Q$ is
$\xso(2,2)_0/\langle-I\rangle=\xpsl_2\xr\times\xpsl_2\xr$. The action
of the left factor preserves each line
$L_{(\lambda:\mu)}=\{(A:B:C:D) \mid (A:B)=(C:D)=(\lambda:\mu)\}$ while the
right action preserves each line
$R_{(\lambda:\mu)}=\{(A:B:C:D) \mid (A:C)=(B:D)=(\lambda:\mu)\}$. We call
the region $AD-BC>0$ \textit{anti de~Sitter space}. It carries a
canonical Lorentzian structure, and we let $X$ denote anti de~Sitter
space with its Lorentzian structure. Given a point $x$ of anti
de~Sitter space, draw the tangent cone to the quadric $Q$. This
determines a Lorentzian quadratic form on the tangent space at $x$, up
to scale: The tangents to $Q$ are in the null directions. In fact the
scale of the metric is also determined. Given a spacelike tangent
vector at $x$, the straight line through $x$ meets $Q$ in two points
$E,H$. Given two points $F,G$ in anti de~Sitter space on the line $l$,
the distance $d(F,G)$ can be defined in terms of the cross-ratio of
$E,F,G,H$ but we will not need the explicit formula.

Letting $F,G$ approach $x$ one defines the length of a tangent vector
at $x$. (The lengths of timelike vectors can be defined by
complexifying so that every line is tangent to $Q$ or else meets $Q$
in two distinct points, possibly imaginary.)

This description of a geometry in terms of projective space with a
distinguished quadric is due to Cayley, in the case of elliptic and
hyperbolic geometry. See ~\cite{37}. See ~\cite{38} for the Lorentzian
case.

We may identify anti de~Sitter space with $\xpsl_2\xr$ by identifying
$(a:b:c:d)$ with its representative, the matrix
$\begin{pmatrix}a&b\\c&d\end{pmatrix}$ unique up to multiplication by
  $-1$, with $ad-bc=1$. Then the left and right actions of
  $\xpsl_2\xr$ are by left and right matrix multiplication. The
  Killing form is bi-invariant, which gives an alternative construction
  of the Lorentzian pseudometric. So we have a model of the Lorentzian
  symmetric space $\xpsl_2\xr$ in which geodesics are represented by
  straight lines in $\xrp^3$. Totally geodesic hyperplanes are the
  intersections of planes with de~Sitter space. So they are given by
  linear equations $P_{(e:f:g:h)}=\{(a:b:c:d):ae+bf+cg+dh=0\}$. A
  plane is spacelike, null or Lorentzian according as
  $eh-gf>0,\;=0,\text{ or }<0$. All the geodesics which are normal to 
$P_{(e:f:g:h)}$ meet in the point $(e:f:g:h)$. In particular if
  $P_{(e:f:g:h)}$ is a null plane, $P_{(e:f:g:h)}$ is tangent to $Q$
  at $(e:f:g:h)$. Each spacelike plane is isometric to the hyperbolic
  plane, each null plane has a unique ruling by null lines and the
  metric is a transverse metric to this foliation, and each Lorentzian
  plane is topologically a M\"obius strip and is isometric to $1+1$-dimensional de~Sitter space. Now let $M$ be a Lorentzian
  manifold locally isometric to anti de~Sitter space and containing a
  closed spacelike hypersurface $S$ with an orthochronous
  neighbourhood. There is a development map $d:\xs\to\xpsl_2\xr$. Given
  a future pointing null direction $\mathbf{u}$ at $x\in\xs$, there is
  a null geodesic at $d(x)$ with tangent vector $d_\ast\mathbf{u}$
  where $d_\ast$ is the derivative of $d$. This null geodesic is
  tangent to $Q$ at a point lying on a line $L_{a(\mathbf{u})}$ of the
  left ruling and a line $R_{b(\mathbf{u})}$ of the right ruling. This
  defines two trivializations of the bundle of future pointing null
  directions. Now the lines of the left, respectively the right
  rulings are permuted by the right, respectively the left action of $\xpsl_2\xr$.

\begin{propo}\label{p19}
The holonomy $\rho=(\rho_L,\rho_R):\pi_1S\to\xpsl_2\xr\times\xpsl_2\xr$
satisfies the condition that $\rho_L,\rho_R$ have Euler class $2-2g$
and therefore are isomorphisms to Fuchsian groups.
\end{propo}

\begin{proof}
As in proposition ~\ref{p1}, the unit tangent bundle may be identified
with the bundle of future pointing null directions. So the maps $a,b$
introduced above identify $UT(S)$ with the flat $\xrp^1$-bundles
associated with $\rho_R:\pi_1S\to\xpsl_2\xr$, respectively
$\rho_L:\pi_1S\to\xpsl_2\xr$. So each of $\rho_L,\rho_R$ has Euler class
$2-2g$ and by Goldman's theorem $\rho_L,\rho_R$ are isomorphisms to
cocompact Fuchsian groups.
\end{proof}

So the holonomy of a spacelike surface $S$ in an anti de~Sitter
spacetime can be considered as a point of $\xteich(S)\times\xteich(S)$. 
(Note the analogy with quasifuchsian groups). Now we will show that 
any point in $\xteich(S)\times\xteich(S)$ may be attained.

For the moment let us distinguish $G_L$, the subgroup of the identity
component $\xaut_0X$ which acts by matrix multiplication on the left,
from $G_R$, the subgroup which acts by matrix multiplication on the
right. We will want to think of $\xaut_0X$ acting on the right of
$X$. Thus the action $G_L\times G_R\times X\to X$ is given by
$(g,h,x)\mapsto g^{-1}\cdot x\cdot h$. Recall that $G_L$ permutes the
lines of the right ruling, which are parametrized by a projective line
$\xrp_{L}^{1}$ so $G_L$ acts on $\xrp_{L}^{1}$ and similarly $G_R$
acts on $\xrp_{R}^{1}$. Now consider the intersection $P\cap Q$ of the
quadric $Q$ with a spacelike plane $P$. $P\cap Q$ is a conic and
determines a one-to-one correspondence between the set of lines of the
left ruling and the set of lines of the right ruling, since $P\cap Q$
meets each line in one point. So $\xrp_{L}^{1}$ and $\xrp_{R}^{1}$ can
be identified with $P\cap Q$. Under this identification $G_L$ and
$G_R$ are each identified with the identity component of the group of
projectivities of the plane $P$ which preserve $P\cap Q$. So $P$
determines an isomorphism $\phi:G_R\to G_L$. The graph of this
automorphism is a subgroup
$\xaut_PX=\{\alpha\in\xaut_0X:\alpha(x)=\phi(g)^{-1}\cdot x\cdot g\}$
of $\xaut_0X$. If we identify $G_L,G_R$ with $\xpsl_2\xr$ using matrix
multiplication, $\phi(g)=aga^{-1}$ for some $a\in\xpsl_2\xr$. If
$\alpha\in \xaut_PX$, $\alpha(x)=ag^{-1}a^{-1}xg$; that is, $a^{-1}x$
is conjugated, so $a$ is fixed by $\alpha$. Thus each $\xaut_PX$ is
the stabilizer of a certain point $a(p)\in X$; $a(P)$ is the dual
point of $P$, i.e., the intersection of all the normal geodesics to
$P\cap X$.

Now given two homomorphisms $\rho_L:\pi_1S\to
G_L=\xaut_0\xrp^{1}_{L},\;\rho_R:\pi_1S\to G_R=\xaut_0\xrp^{1}_{R}$,
each with Euler class $2-2g$, there is a homeomorphism
$h:\xrp^{1}_{L}\to\xrp^{1}_{R}$ conjugating $\rho_L$ to $\rho_R$, because
$\rho_L(\pi_1S),\rho_{R}(\pi_1S)$ are isomorphic cocompact Fuchsian
groups and so their actions on their limit circles are topologically
conjugate. Moreover, if $\xrp_{L}^{1}$ and $\xrp_{R}^{1}$ are oriented
so that the map $i_p:\xrp_{L}^{1}\to\xrp^{1}_{R}$ determined by a
spacelike plane $P$ is orientation preserving, the homeomorphism $h$
is orientation preserving. We may consider the graph of $h$ as
embedded in $Q$, as $Q\cong\xrp_{L}^{1}\times\xrp_{R}^{1}$. $Q$ has a
conformally Lorentzian structure. Because $h$ is orientation
preserving, nearby points on the graph of $h$ are spacelike with
respect to each other. (For an orientation reversing homeomorphism,
nearby points would be timelike with respect to each other. The
intersection of $Q$ with a timelike plane is the graph of an element
of the non-identity component of $\xpgl_2\xr$, if a spacelike plane has
been used to identify $\xpsl_2\xr_L$ and $\xpsl_2\xr_R$.) So given 3
points $A,B,C$ on $\xgraph(h)$, the plane through $A,B,C$ is
spacelike.

\begin{lemma}\label{l5}
There is a plane (necessarily spacelike) disjoint from
$\xgraph(h)\subset Q$.
\end{lemma}

\begin{proof}
Given 3 points $A,B,C$ on the graph of $h$ let $A',B',C'$ be the
intersections of lines $l_A,l_B,l_C$ of the left ruling with lines 
$r_B,r_C,r_A$ respectively of the right ruling, and let $A'',B'',C''$
be the intersections of $l_A,l_B,l_C$ with $r_C,r_A,r_B$. Then the
graph of $h$ lies in the three twisted quadrilaterals
$AA'BB'',BB'CC'',CC'AA''$ on $Q$.  The configuration is unique up to a projective transformation and can be realized in $\xr^3$ with $Q$ given by 
$x^2+y^2=z^2+1$ and
\begin{align*}
A&=(1,0,0), B=(\cos\frac{2\pi}{3},\sin\frac{2\pi}{3},0), C=(\cos\frac{2\pi}{3},-\sin\frac{2\pi}{3},0),\\
A'&=(2\sin\frac{\pi}{3},2\cos\frac{\pi}{3},\sqrt{3}), B'=(-2,0,\sqrt{3}), C'=(-2\sin\frac{\pi}{3},2\cos\frac{\pi}{3},\sqrt{3})
\end{align*}
and $A'',B'',C''$ are the reflections of $C',A',B'$ in the plane
$z=0$. Then there are many spacelike planes, including the plane at
infinity and the planes $z=k$ where $|k|>\sqrt{3}$, disjoint from $\xgraph(\phi)$.
\end{proof}

Having chosen a spacelike plane $P_\infty$, the convex hull $X(\phi)$
of $\xgraph(\phi)$ in $\xrp^3-P_\infty$ is defined, and $\pi_1S$ acts
on $X(\phi)$ by $\rho=(\rho_L,\rho_R):\pi_1S\to G_L\times G_R$. From
the fact that $\xgraph(h)$ is a spacelike curve in the conformally
Lorentzian structure on $Q$ it follows that lines and planes in the
boundary of the convex hull are all spacelike.

\begin{propo}\label{p20}
Given a point $(x,y)$ in $\xteich(S)\times\xteich(S)$ there exists a
compact spacetime $X(x,y)$ with boundary, homeomorphic to
$S\times[0,1]$ by a homeomorphism taking any $S\times\{x\}$ to a
closed spacelike hypersurface, which is locally modelled on anti
de~Sitter space, and which has left and right holonomies $\rho_L,\rho_R$
determined by the point $(x,y)$. Moreover the boundary is spacelike,
locally convex and has no extreme points and these conditions uniquely
determine $X(x,y)$ up to Lorentz isometry.
\end{propo}

\begin{proof}

Choose homomorphisms $\rho_L,\rho_R:\pi_1S\to\xpsl_2\xr$ representing
$(x,y)\in\xteich(S)\times\xteich(S)$. We have seen that there exists a
compact convex set $X(\phi)$, equivariant with respect to
$(\rho_L,\rho_R):\pi_1S\to G_L\times G_R=\xaut_0X$. The boundary of
$X(\phi)$ in $Q$ is the graph of a homeomorphism
conjugating $\rho_L$ to $\rho_R$, and $F(\phi),P(\phi)$ are the
future and past boundary components of $X(\phi)$. $F(\phi),P(\phi)$
are locally convex and spacelike.

We need to see that $\rho(\pi_1S)$ acts properly discontinuously on
$F(\phi)$. Here is one way to see this; proposition ~\ref{p21} will
give another. Given $p\in F(\phi)$ and a sequence of distinct group
elements $\gamma_i\in\pi_1S$, after passing to a subsequence, there
exists $x\in\xrp_{R}^{1}$ such that $\rho_R(\gamma_i y)\to x\in
\xh^2\cup\xrp_{L}^{1}$, uniformly for $y$ in a compact subset of
$\xh^2$. Fix a spacelike supporting plane $T$ through $p$. (If the
supporting plane is unique, it is the tangent plane.) There is a
rotationally invariant measure on the rays from $T$, and this
determines a measure $\mu_{p,T}^{L}$ on $\xrp_{L}^{1}$ and a measure 
$\mu_{p,T}^{R}$ on $\xrp_{R}^{1}$. Now
$\rho_L(\gamma_i)_\ast\mu_{p,T}^{L}\to\delta_{\phi(x)}\in
M(\xrp_{L}^{1})$ and $\rho_R(\gamma_i)_\ast\mu_{p,T}^{R}\to\delta_{x}\in
M(\xrp_{R}^{1})$ where $M(\xrp_{L}^{1}),M(\xrp_{R}^{1})$ are the
spaces of probability measures on $\xrp_{L}^{1},\xrp_{R}^{1}$
respectively and where $\delta_{\phi(x)},\delta_x$ are the point
masses at $\phi(x)$ and $x$ respectively. Now suppose that
$\rho(\pi_1S)$ did not act properly discontinuously on $F(\phi)$. Then
for some $p$, there is a sequence $\gamma_i$ such that
$\rho(\gamma_i)\cdot p$ remains within a compact set $K\subset
F(\phi)$. But then the measures $\rho_{L}(\gamma_i)_\ast\mu_{p,T}^{L}$
lie in a bounded set in the space $C^0(\xrp_{L}^{1})$, the space of
measures with a continuous density with respect to the
$\xso(2)$-invariant measure, and similarly for
$\rho_{R}(\gamma_i)_\ast\mu_{p,T}^{R}$. So $\rho(\pi_1S)$ acts
properly discontinuously on $F(\phi)$.

To deduce that the action on $X(\phi)$ is properly discontinuous,
observe that if $p\in X(\phi)$, the future pointing null and timelike
rays from $p$ meet $F(\phi)$ in a compact disc $D(p)$. $D(p)$ depends
continuously on $p$.

Now for the uniqueness of $X(\phi)$.

\begin{lemma}\label{l6}
Suppose $f:A\to X$ is a locally isometric map of a complete connected
Riemannian manifold $A$ to anti de~Sitter space. Then a) $f$ is
proper; b) $f$ is an embedding, and $f(A)$ intersects every timelike
geodesic in only one point.
\end{lemma}

\begin{proof}
Fix an action of $\xso(2)$ on $X$ by right multiplication, and
consider the Riemannian submersion $p:X\to X/\xso(2)=\xh^2$. (The
preimage by $p$ of a round disc in the hyperbolic plane is a solid
torus with boundary a doubly ruled hyperboloid; one of the rulings is
by orbits of the $\xso(2)$-action.) The derivative of $p$ can only
lengthen a spacelike vector. (Note however that $X$ does not contain a
subset which is mapped isometrically to $\xh^2$ by $p$, because the
field of planes orthogonal to the $\xso(2)$-action is totally
non-integrable. In particular the restriction of $p$ to a totally
geodesic spacelike hypersurface, which is isometric to $\xh^2$ is not
an isometry.) So if $A$ is complete, the projection must be a covering
map. Since $\xh^2$ is simply connected, $A$ is diffeomorphic to
$\xh^2$, $f$ is proper. Given any timelike geodesic, there is a unique
$\xso(2)$-action by right multiplication of which it is an orbit, so
$b)$ follows.
\end{proof}

Now we identify $A$ with its image, and consider the closure $\xcl A$
in $\xrp^3$ of $f(A)$. We want to show $\xcl A$ is a disc.

In lemma ~\ref{l7}, we consider that a null geodesic $l\cap X$ does
not contain its point at infinity, the point of tangency $l\cap Q$ of
tangency of the line $l$ with $Q$. Thus a point $x$ divides a null
geodesic into two half lines, the past and future pointing null rays
through $x$, each of which is homeomorphic to $[0,\infty]$.

\begin{lemma}\label{l7}
a) Given a spacelike disc $D$ which meets no timelike or null geodesic
in more than one point, let $F(D)$, respectively $P(D)$, be the set of
$x\in X$ such that every past pointing, respectively future pointing
null ray or timelike geodesic through $x$ meets $D$. Then $P(D)\cap
F(D)=D$. b) Given a complete spacelike submanifold $A\subset X$, the
closure of $A$ in $X\cup Q\subset\xrp^3$ is a disc $A\cup\partial A$
and the boundary $\partial A=\xcl A\cap Q$ is a circle which is
nowhere timelike in the conformally Lorentzian structure on $Q$.
\end{lemma}

\begin{figure}[t!]
\begin{center}
\epsfysize=3in
\epsffile{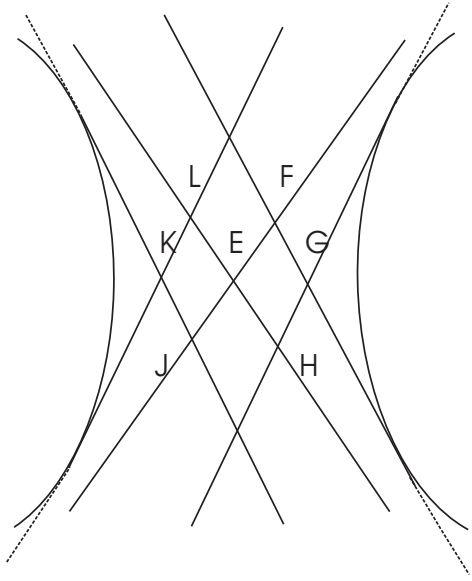}
\end{center}

\caption{\label{fig5}A nowhere timelike curve on the quadric at infinity}
\end{figure}

Before proving the lemma, let us amplify, perhaps pedantically, the
statement of part $b)$. A differentiable curve on $Q$ is nowhere
timelike if its tangent vector is nowhere timelike in the conformally
Lorentzian structure on $Q$. In general, a nowhere timelike curve is a
curve that can be approximated by differentiable nowhere timelike
curves. An alternative definition may be given without recourse to
differentiability: A curve $C:t\to C(t)$ on $Q$ is nowhere timelike if
given any point $E=C(s)$ on $C$ and any sufficiently small pair
$EFGH,EJKL$ if quadrangles on $Q$ such that $JEF,HG,KL$ are short
segments of lines of the left ruling and $GF,HEJ,JK$ are short
segments of lines of the right ruling, $C(t)$ lies in the union of the
two small quadrangles $EFGH,EJKL$ for any $t$ sufficiently close to
$s$.  See Figure~\ref{fig5}. Note that given any two points on $Q$
which do not lie on a line in one of the two rulings, there is a
spacelike geodesic in $X$ with the two given points as endpoints.

\begin{proof}

$a)$ If there was a point $y\in P(D)\cap F(D)$ such that $y$ is not in
  $D$ then any null geodesic through $y$ meets $D$ in two points. This
  proves $a)$. Now let $H(D)=P(D)\cup F(D)$.

\begin{sublemma}\label{sl1}
Given any complete spacelike surface $A\subset X$, there is a
spacelike plane in $X$ whose closure in $X\cup Q \subset \xrp^3$ is disjoint
from the closure of $A$.
\end{sublemma}

\begin{proof}
Given $p\in A$ consider the dual plane to $p$, which meets $Q$ in the
points of tangency of the null lines through $p$. Every normal to the
dual plane passes through $p$, so if $A$ met the dual plane there
would be a timelike geodesic which met $A$ in two distinct points. To
see that the dual plane is disjoint from the closure of $A$ in
$\xrp^3$, assume without loss of generality that $A$ contains a
neighbourhood of $p$ which lies in a plane. The union of the dual
planes is a neighbourhood of the dual plane of $p$.
\end{proof}

Thus we can choose a plane at infinity, and then the convex hull of
the closure of $A$ in $\xrp^3$ is compact. Furthermore, there is a
compact convex set which is an intersection of half-spaces with
spacelike boundaries, and contains $A$; let $B=B(A)$ be the minimal
such set. A half-space with spacelike boundary can either contain the
immediate past of its boundary in which case we will call it
``past-complete'' or not; the boundary of $B$ divides into the future
boundary where every point has at least one spacelike supporting plane
which $B$ lies in the past of, the past boundary, and the rest.

Fix a point $x\in X$. Let $D_r$ be the disc in
$\xh^2=X/\xso(2)$ with radius $r$ and center $p(x)$, and let
$T_r=p^{-1}D_r$ be its preimage in $X$ and $A_r=A\cap T_r$. There is
an $\xso(2)\times\xso(2)$ action on $X\cup Q$ whose orbits are the
tori $T_r$. A spacelike line through $x$ together with this
$\xso(2)\times\xso(2)$ action determines identifications of the tori
$T_r$ with $Q$. Choose a sequence $r_n\to\infty$ such that $\partial
A_{r_n}$ converges in the space of compact subsets of $\xrp^3$ with
the Hausdorff topology, or equivalently, in the Hausdorff topology on
the space of compact subsets of $Q$ using the identifications of the
tori $T_r$ with $Q$. Each $\partial A_{r_n}$ can be considered as the
graph of a homeomorphism between the left and right copies of
$\xso(2)$. So the relation $\partial A_\infty$ can be approximated by
homeomorphisms. It follows that $\partial A_\infty$ intersects each
line on $Q$ in a nonempty connected set. It follows that $\partial
A_\infty$ is a nowhere timelike topological circle. There is a plane
disjoint from $\partial A_\infty$ by (the argument of) lemma
~\ref{l5}. Consider the sets $H(A_r)$. Since this is an increasing
family, parametrized by $[0,\infty]$, of compact sets any subsequence
has the same Hausdorff limit. Let us define $H(A_\infty)$ to be the
set of points which have dual planes which do not cross the circle
$A_\infty$. (The circle $A_\infty$ separates a small neighbourhood of
itself in $Q$. Let $N(y)$ be the intersection of the dual plane of $y$
with $Q$. $N(y)$ crosses $A_\infty$ at $w\in A_\infty$ if $N(y)$
contains arbitrarily small intervals containing $w$ which meet both
components of the complement of $A_\infty$ in the small neighbourhood
of $A_\infty$.) Then the closed set $H(A_\infty)$ contains all the
sets $H(A_r)$. For if $y\in  H(A_r)$ the null cone of $y$ intersects
$A$ in a compact set, and so the conic say $N(y)$ of tangency of the
null cone with $Q$ can't cross $A_\infty$ because if $N(y)$ did then
the null cone would intersect $A_{r_n}$ for all sufficiently large
$n$. Since $H(A_\infty\cap Q)=A_\infty$, the Hausdorff limit of the
sets $H(A_r)$ is $H(A_\infty)=\partial A$ and also the closure of $A$
in $X\cap Q$ is $A\cap\partial A$. Now to see that $\xcl
A=A\cap\partial A$ is a disc, consider the action of the diagonal
subgroup of $\xso(2)\times\xso(2)$ on $X\cap Q$. It gives a continuous
bijection of $\xcl A$ onto the disc. Since $\xcl A$ is compact this is
a homeomorphism.
\end{proof}

The following sublemma, the proof of which follows a suggestion of
B.~Bowditch and N.H.~Kuiper, can be used to simplify the proof of
lemma ~\ref{l7}. It is of some independent interest; perhaps it is well-known.

\begin{sublemma}\label{sl2}
Given an immersion $i:D^n\to\xrp^n$ of a disc such that the boundary
is strictly locally convex, $i$ is an embedding and the image is a
convex set in the complement of some hyperplane.
\end{sublemma}

\begin{proof}
It follows from the local convexity that any two points can be joined
by a distance minimizing geodesic. $i$ is an embedding because if
$p,q\in D^n$ and $i(p)=i(q)$, then for any $r$ in $\partial D^n$ such
that $i(r)\ne i(p)$ the geodesics from $r$ to $p,q$ must have a common
initial segment, which is in the interior of $D^n$ by local
convexity. So one (say $p$) of $p,q$ is in the interior. As every
geodesic ray through $p$ is the preimage of a closed geodesic, every
geodesic ray through $p$ eventually reaches the boundary or else
covers a closed geodesic. Since every ray through $p$ which reaches
the boundary passes through $q$, possibly after extension through $p$,
either there is a closed geodesic or else an open set of geodesic rays
through $p$ which reach $q$. The set of geodesics from $p$ which reach
$q$ must be open as well as closed because of the convexity of the
boundary. It follows that all geodesics from $p$ reach $q$ and then
return to $p$, so $D^n$ has no boundary, which is absurd. Now consider
$D^n$ as embedded in the double cover of $\xrp^n$ and consider the
cone over $D^n$ in $\xr^{n+1}$; it is strictly convex and it follows
that $D^n$ lies on one side of a hyperplane and is a convex set in the
affine space defined as the complement of that hyperplane.
\end{proof}

Here is an alternative proof of lemma ~\ref{l7}, worked out with
Bowditch. $\partial T_r$ is an orbit of $\xso(2)\times\xso(2)$ and 
therefore is a doubly ruled quadric. So the interior of $T_r$ has a 
Lorentzian structure defined by $\partial T_r$. In this Lorentzian
structure, 
light cones are narrower than in the Lorentzian structure of $X$, so 
more vectors are spacelike. So the intersection $A_r=A\cap T_r$ is
spacelike 
in the new Lorentzian structure. After choosing a different product
structure 
on $\xso(2)\times\xso(2)$ the intersections $A_r$ become graphs of
functions 
all of which have Lipschitz constant $1$. So it suffices to show
pointwise 
convergence at one point in order to show that the sets $A_r$ converge 
to a nowhere timelike circle. It is easy to see that the intersection
of $A$ 
with any indefinite plane has closure an arc, and the pointwise
convergence 
follows. This proof is shorter, but it seems useful to understand the
sets 
$H(A_r)$ occurring in lemma ~\ref{l7}.

Extending the definitions of $F(D),P(D),H(D)$ given in lemma ~\ref{l7} 
and its proof to the case of a noncompact spacelike surface, we see
that 
$H(A)=H(A_\infty)$ where $H(A_\infty)$ has just been defined. 
Suppose $y$ lies in the interior of the convex hull of $\partial
A$. Then no point in the dual plane of $y$ can be in $H(A)$. For if
$z$ is in the dual 
plane of $y$, then the dual plane of $z$ passes through $y$, and so
meets 
$Q$ in a conic which crosses $\partial A$. Thus there is a plane
disjoint 
from $H(A)$, except in the case that $A$ is a spacelike plane.  In any case there is a plane disjoint from the interior of $H(A)$.  Supposing there is a plane disjoint from $H(A)$, $H(A)$ is the intersection of all the closed half-spaces
with boundary a totally geodesic null plane which contain $A$. So
$H(A)$ is a compact convex set, except when $A$ is a plane in which
case $H(A)$ is 
the union of a closed convex cylinder and one point at infinity. 
We can also think of $H(A)$ as the closure of the region in which 
any spacelike surface with boundary $\partial A$ must lie. 
Now given a compact Lorentzian manifold $Z$, locally isometric 
to anti de~Sitter space, with spacelike locally convex boundary with 
no extreme points, with the holonomy of $X(\phi)$, the manifold can be 
slightly thickened to $Z'$ to make the boundary smooth. 
Then lemma ~\ref{l7} applies, and we deduce that the closure of the 
universal cover of $\partial Z'$ is the union of $\xgraph(\phi)$ and
two 
locally convex spacelike surfaces. 
Alternatively lemma ~\ref{l7} could have been proven for spacelike 
surfaces which were not smoothly immersed. 
It follows that the universal cover of either boundary component of
$Z$ 
is one of the two boundary components of $X(\phi)$. 
So the boundary is $F(\phi)\cup P(\phi)$. 
It follows that $Z$ is the quotient of $X(\phi)$.
\end{proof}

\begin{propo}\label{p21}
Given an orthochronous anti de~Sitter spacetime which contains a 
product neighbourhood $N(F)$ of a closed spacelike surface $F$, 
then there is a smaller neighbourhood $N'(F)$ which embeds in 
$\xint\big(H(F(\phi))\big)/\pi_1F$ where $\phi$ conjugates the left 
and right holonomies of $\pi_1F$. 
If $N(F)$ is a domain of dependence then $N(F)$ embeds in 
$\xint\big(H(F(\phi))\big)/\pi_1F$. 
So $\xint\big(H(F(\phi))\big)/\pi_1F$ is a domain of dependence 
which contains all other domains of dependence with the same holonomy.
\end{propo}

\begin{proof}
First we show that $F$ has negative Euler characteristic. 
$F$ cannot be a sphere or projective plane, because anti 
de~Sitter space contains no spacelike sphere. 
If $F$ had Euler characteristic $0$, then $F$ or a double cover is a
torus; 
suppose $F$ is a torus. By lemma ~\ref{l6} the holonomy cover $\widehat{F}$ of 
$F$ is embedded and simply connected. 
We could argue that $\widehat{F}$ is
quasi-isometric 
to a hyperbolic plane and so the ball of radius $r$ in $\widehat{F}$ 
grows exponentially with $r$, while the cover of a torus has only 
polynomial volume growth. Alternatively, the left holonomy defines a 
homomorphism from $\xz\oplus\xz$ into an abelian subgroup, 
necessarily a $1$-parameter subgroup, of $G_L$. 
So after an arbitrarily small change in holonomy the image is cyclic. 
The right holonomy is topologically conjugate to the left holonomy so
the holonomy has kernel which contradicts lemma ~\ref{l6}. 
So $F$ has negative Euler characteristic.

We know that there is a manifold $X(\phi)/\pi_1F$ with the same
holonomy as $N(F)$, and it is easy to see that $\pi_1F$ acts properly 
discontinuously on the enlargement
$\xint\big(H(F(\phi))\big)=\xint\big(H(P(\phi))\big)=\xint(H(\xgraph\phi))$. 
The universal cover of $F$ isometrically immerses in anti de~Sitter
space, 
equivariantly with respect to the action of the holonomy of
$X(\phi)/\pi_1F$, 
and we know that this immersion is an embedding and the closure of the 
embedding is a disc. The boundary of the disc can be thought of as a 
relation conjugating the actions of $\pi_1F$ on the parameter space 
$\xrp_{L}^{1}$ and $\xrp_{R}^{1}$ of the left and right rulings. 
Such a relation must be a homeomorphism (so the boundary of the discs 
contains no segment of a line of $Q$) and $\phi$ is 
the unique homeomorphism which conjugates the left and right actions. 
So the boundary of the closure of the universal cover of $F$ is the 
boundary $\xgraph\phi$ of $X(\phi)$. 
Therefore the universal cover of $F$ is equivariantly embedded in 
$\xint\big(H(F(\phi))\big)=\xint(H(\xgraph\phi))$. 
So some neighbourhood of $F$ embeds in
$\xint\big(H(F(\phi))\big)/\pi_1F$. 
If $N(F)$ is a domain of dependence then the development map must 
define a map taking $N(F)$ to $\xint\big(H(F(\phi))\big)/\pi_1F$. 
If this were not injective there would be two points in the 
universal cover of the domain of dependence whose timelike and null 
geodesics hit the same subset of $\widetilde{F}$, but this is impossible.
\end{proof}

We can reinterpret the boundary components of the convex hull as
earthquakes, ~\cite{14,40,41,20} of which we will now give a fairly
self-contained exposition. The simplest earthquakes are fractional
Dehn twists, i.e., shears on closed geodesics; a geodesic lamination
is a limiting case of a long simple closed geodesic and an earthquake
is a shearing motion along a geodesic lamination. Recall that
$1+1$-dimensional de~Sitter space $\xrp^2-\xcl\xh^2$ can be identified
with the set of (unoriented) geodesics: Given $p\in\xrp^2-\xcl\xh^2$
draw the two tangents to the conic $S^{1}_{\infty}$, and identify $p$
with the geodesic joining the two points of tangency. The geodesics
corresponding to $p$ and $q$ cross, respectively are asymptotic, iff
the line joining $p$ and $q$ does not meet  $S^{1}_{\infty}$,
respectively is tangent to  $S^{1}_{\infty}$, iff $p$ and $q$ are
mutually spacelike, respectively mutually null. A
geodesic lamination is a partition of a closed subset of $\xh^2$ into
(disjoint) complete geodesics. Dually it is a closed subset of
$\xrp^2-\xcl\xh^2$ containing no mutually spacelike pair of points. A
measured lamination $(\lambda,\mu)$ is a closed subset $\lambda$ of
$\xrp^2-\xcl\xh^2$ together with a locally finite positive measure
$\mu$ such that $\lambda=\xsupp\mu$. Dually it is a geodesic
lamination together with a measure defined on transverse arcs. Given a
measured lamination $(\lambda,\mu)$, let $\lambda_0$ be the union of
the leaves of $\lambda$ which carry atoms of the measure
$\mu$. $\lambda_0$ consists of countably many leaves, but need not be
closed, though if in addition $\lambda$ is invariant under a Fuchsian
group with finite covolume, $\lambda_0$ will be closed
~\cite{19}. There is a metric space $\xh^2(\lambda_0)$ in which
distances are given by lengths of paths together with a map
$p_{\lambda_0}:\xh^2(\lambda_0)\to\xh^2$ such that if $l\in\lambda_0$
is a geodesic then $p_{\lambda_0}^{-1}(l)$ is isometric to
$l\times[0,\mu(l)]$ and $p_{\lambda_0}$ is the projection onto the
factor $l$. $\xh^2(\lambda_0)$ is unique up to a unique isometry, and
homeomorphic to $\xr^2$. $\xh^2(\lambda_0)$ is
complete. $\xh^2(\lambda_0)$ carries a measured geodesic lamination
$\lambda^\ast$ without atoms (generalizing the definition in an
obvious way) and the push forward of the transverse measure of
$\lambda^\ast$ equals the transverse measure $\lambda$.

A \textit{left earthquake} with \textit{shearing lamination} $\lambda$
is a relation $f:\xh^2\to\xh^2$ such that 

$i)$ $f=f^\ast\circ
p_{\lambda_0}^{-1}$ where $f^\ast:\xh^2(\lambda_0)\to\xh^2$ is a
function,

$ii)$ $f^\ast|_{p_{\lambda_0}^{-1}(l)}$ maps $p_{\lambda_0}^{-1}(l)$
to a geodesic $f(l)$. $f^\ast$ is an isometry on each $l\times x$ and
each $q\times[0,\mu(l)]$. Given a transverse arc
$a:(-\epsilon,\epsilon)\to\xh^2$ parametrized by arc length, with
$a(0)\in l$ for some geodesic $l$ of the lamination $\lambda$, there
is an obvious lift $a^\ast:\big(-\epsilon,\epsilon+\mu(\lambda_0\cap
a(-\epsilon,\epsilon))\big)\to \xh^2(\lambda_0)$. It is required that
if $l$ is oriented so as to cross $a((-\epsilon,\epsilon))$ from right
to left looking in the positive direction along
$a((-\epsilon,\epsilon))$, 
then the
lift $a^\ast|_{[0,\mu(l)]}$ be orientation preserving: $a^\ast(t)$
moves to the left, as seen looking in the positive direction along
$a((-\epsilon,\epsilon))$, as $t$ increases; and

$iii)$ $f^\ast(x)=x\cdot T(x)$, where
$T:\xh^2(\lambda_0)\to\xiso^+(\xh^2)$ is continuous and, if
$a:[0,\epsilon)\to\xh^2(\lambda_0)$ is a transverse arc,
  $T(a(x))=h(\mu^\ast(a[0,x]))$ where $\frac{d}{dt}|_{t=0}\frac{h(t)}{h(0)}=t_l$
  where $t_l$ is the element of the Lie algebra of $\xiso^+(\xh^2)$
  which translates $l$ along $l$ with unit speed toward the left, as
  seen looking in the positive direction along
  $a((-\epsilon,\epsilon))$.

Now let $\phi:S^1\to S^1$ be a homeomorphism. Embed $S^1$ in $Q$ as
the intersection $C$ of $Q$ with a spacelike plane $P$. We identify
$q\in Q$ with $(a,b)\in\xrp_{L}^{1}\times\xrp_{R}^{1}$ where $a,b$ are
the lines through $q$ belonging to the left and right rulings. We
identify each of $\xrp_{L}^{1},\xrp_{R}^{1}$ with $C$ by identifying
$a\in\xrp_{L}^{1}$ or $\xrp_{R}^{1}$ with $a\cap C$. Then
$\xgraph(\phi)=\{a\in Q: a=(x,\phi(x))\text{ for some }x\}$. Now
consider the future component $F(\phi)$ of the boundary of the convex
hull of $\xgraph(\phi)$. It carries a lamination $\lambda'$ (i.e., a
closed subset partitioned into properly embedded lines): The union of
those spacelike geodesics in anti de~Sitter space which lie in
$F(\phi)$ but do not lie in the interior of the intersection of a
spacelike plane with $F(\phi)$. For each geodesic $l'$ in this
lamination, with end points $a',b'$ let $a,b$ be the intersections of
the geodesic in the left ruling through $a,b$ respectively with the
conic $C$. Let $l$ be the geodesic through $a,b$. The union of the
geodesics $l$ is a geodesic lamination $\lambda$ on $P$. Similarly
using the right ruling we construct a geodesic lamination $\lambda''$ on $P$. Now
given a supporting plane $T$ of $F(\phi)$ there is a unique
element $g(T)$ of $G_L$ such that $T\cdot g(T)=P$. $g(T)$ moves each
point of $T\cap Q$ along a line of the left ruling to the intersection
of that line with $P$. $g(T)$ is also the unique element of $G_L$
which translates the dual point of $T$ to the dual point of
$P$. Similarly there is a unique element $h(T)\in G_R$. Now we define
a left earthquake, which extends continuously to the boundary of
$\xh^2$ where it equals $\phi$. Given a supporting plane $T$ of
$F(\phi)$, map $T\cap F(\phi)\cdot g(T)$ to $T\cap F(\phi)\cdot h(T)$
by $(g(T)^{-1},h(T))$. This defines a relation which is single valued
except on sets $T\cap F(\phi)\cdot g(T)$ for which $T$ is not the
unique supporting plane at $T\cap F(\phi)$. The lamination $\lambda'$
carries a bending measure (discussed at length in ~\cite{44,20}; there
is no essential difference in the case of an indefinite metric) which
can be transported to the laminations $\lambda,\lambda''$. One way to
describe the bending measure is to associate to a transverse arc
the length of the arc of dual points to the supporting
planes. To do so one needs to know that the arc of dual points is
rectifiable. This follows however from the fact that the direction of
the geodesic lamination $\lambda$ is a Lipschitz continuous 
function ~\cite{19}.

\begin{propo}\label{p22}
For each orientation preserving homeomorphism $\phi:S^1\to S^1$ there
is a unique left earthquake $E(\phi):\xh^2\to\xh^2$ which extends
continuously to $\xh^2\cap S^1$ with value $\phi$ on $S^1$. The map
$\phi\mapsto E(\phi)$ is equivariant with respect to pre-composition
and post-composition with $\xiso^+\xh^2$. The shearing lamination and
measure are naturally associated to the bending lamination and bending measure of the boundary of the convex hull in anti de~Sitter space of
the graph of $\phi$ regarded as a subset of $Q$. There is also a
unique right earthquake obtained from the past boundary of the convex hull.
\end{propo}

\begin{figure}[t!]

\begin{center}
\epsfysize=3in
\epsffile{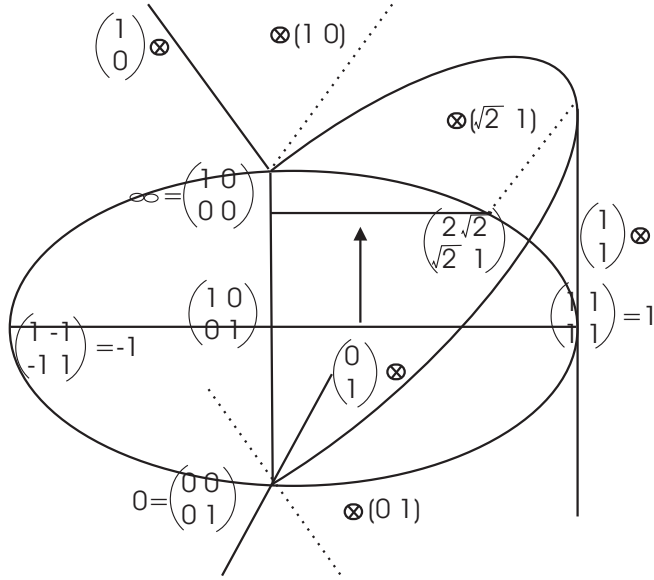}
\end{center}

\caption{\label{fig6}A left earthquake}
\end{figure}

\begin{proof}
One need only check the definitions, except perhaps for the
uniqueness. Given an earthquake, one can construct a locally convex
surface with no extreme points with boundary $\xgraph\phi$; this must
be the future boundary. The easiest way to check that the future and
past boundaries correspond to left and right earthquakes is to draw
some examples. (In Figure ~\ref{fig6}, the earthquake is certainly a
left earthquake; the corresponding component of the boundary of the
convex hull, which is the union of two spacelike half-spaces, may look
at first like the past boundary, but the future direction along say
the subgroup of rotations $\begin{pmatrix}\cos\theta&\sin\theta\\-\sin\theta&\cos\theta\end{pmatrix}$ is downwards
in the picture.)
\end{proof}

Choose $P$ to be the plane $b=c$ in the space of matrices. The
intersection of the plane $P$ with $Q$ can be identified with $\xrp^1$
by $(a:b:c:d)\mapsto(a:b)$. There is an earthquake $E:\xh^2\to\xh^2$
with shearing lamination given by the single geodesic $b=c=0$ and
transverse measure $\log s$. The earthquake can be normalized to be
the identity on the half-plane $ab<0$ and given by 
$A\mapsto\begin{pmatrix}s&0\\0&1\end{pmatrix}A\begin{pmatrix}s&0\\0&1\end{pmatrix}$
on the half-plane $ab>0$. Thus
$\begin{pmatrix}1&1\\1&1\end{pmatrix}\mapsto\begin{pmatrix}s^2&s\\s&1\end{pmatrix}$.
  (For convenience we are using homogeneous coordinates, and matrices
  in $\xgl(2,\xr)$ rather than in $\xsl(2,\xr)$). The boundary value
  of the earthquake is the homeomorphism of $\xrp^1$ which is the
  identity on negative reals and multiplication by $s$ on positive
  reals. The length of the shear along the lamination is $\log s$. 
$\begin{pmatrix}1&1\\1&1\end{pmatrix}$ goes along the line 
$\begin{pmatrix}1\\1\end{pmatrix}\otimes\xr^2$ to
      $\begin{pmatrix}s&1\\s&1\end{pmatrix}$ and then along the line
	$\xr^2\otimes\begin{pmatrix}s&1\end{pmatrix}$ to 
$\begin{pmatrix}s^2&s\\s&1\end{pmatrix}
=E(\begin{pmatrix}1&1\\1&1\end{pmatrix})$. The isometry fixing $b=c=0$
  pointwise and taking the plane through
  $\begin{pmatrix}1&0\\0&1\end{pmatrix},\begin{pmatrix}0&0\\0&1\end{pmatrix},\begin{pmatrix}s&1\\s&1\end{pmatrix}$ to the plane $P$ is given by conjugation by $\begin{pmatrix}s^{\frac{1}{2}}&0\\0&1\end{pmatrix}$ and takes $\begin{pmatrix}s&1\\s&1\end{pmatrix}$ to $\begin{pmatrix}s&s^{\frac{1}{2}}\\s^{\frac{1}{2}}&1\end{pmatrix}$. Thus after an earthquake with half the transverse measure, there is an isometry with the boundary of the convex hull. As long as the earthquake is along a lamination which is a union of isolated leaves, the same argument is applicable: For after a composition with an isometry one can assume the earthquake is the identity on a region on one side of any given leaf. Any complete spacelike surface in anti de~Sitter space which is a boundary component of the convex hull of a nowhere timelike curve in the boundary can be approximated by one with a discrete bending lamination, so in general given a hyperbolic surface $S$ and an earthquake on $S$ determined by a measured geodesic lamination, the quotient of the top of the convex hull in de~Sitter space of the graph on the quadric $Q$ of the conjugating homeomorphism determined by the earthquake on $S$ is isometric to the hyperbolic surface obtained by an earthquake with half as much transverse measure.

\begin{propo}\label{p23}

Suppose $E:\xh^2\cup\xrp^1\to\xh^2\cup\xrp^1$ is a left earthquake
along a lamination $\lambda$ with transverse measure $\mu$. Then there
is an isometry $I$ from the hyperbolic plane to the future component
of the convex hull of the graph of the boundary value homeomorphism
$\phi$ of $E$. Let $F$ be the earthquake with shearing lamination
$\lambda$ and transverse measure $\frac{1}{2}\mu$. Then $I\circ F$ has as
boundary value the map which sends $p\in\partial\xh^2$ to the
corresponding point on the graph of $\phi$, i.e., the intersection of
the line in the left ruling through $p$ with the graph.
\end{propo}

It would be more natural to formulate proposition ~\ref{p23} in terms
of relative hyperbolic structures as in ~\cite{41}. When an
earthquake's shearing lamination (together with its measure) is
equivariant with respect to a Fuchsian group $\pi$ we can consider the
induced earthquake on the quotient surface or orbifold; its image lies
in the quotient surface or orbifold of a different Fuchsian group
$\pi'$. Given a homeomorphism of a non-oriented surface, there is a
canonical lift to the orientation cover: A point $q$ in the
orientation cover is a point $p$ in the non-oriented surface together
with a local orientation at $p$, and the lift maps the $q$ to the
image of $p$ together with the image of the local orientation. Thus
left and right earthquakes can be defined on nonorientable as well as
on orientable surfaces.

Given any measured geodesic lamination $(\lambda',2\mu)$ one can bend
a spacelike geodesic plane along it, obtaining a complete convex
surface $A$ and the closure of $A$ will be a disc in $X\cup Q$ with
nowhere spacelike boundary. (Here anti de~Sitter space has an
advantage over hyperbolic space: No matter how large the bending, the
map to anti de~Sitter space is an isometry.) We regard $\lambda'$ as
lying in $A$. For each geodesic $l'$ in $\lambda'$, there are
corresponding geodesics $l$ and $l''$ in $\xh^2=X\cap P$ obtained by
left and respectively right translation of the endpoints. It follows,
using the previous proposition, that there is an earthquake $E$
defined on $\xh^2$ with shearing lamination $(\lambda',\mu)$; up to an
isometry the image of $\lambda'$ is the ``lamination'' $\lambda''$
which is the union of the geodesics $l''$. $\lambda''$ need not be
closed. The image of $E$ is a convex open set whose frontier in
$\xh^2$ is a union of geodesics. Indeed, not all earthquakes have
boundary values which are homeomorphisms. The graph in $Q$ of
$f:\xr\to\xr,t\mapsto e^t$ union a single segment of a null line has a
convex hull, whose boundary components give earthquakes taking a
complete hyperbolic surface with cyclic fundamental group and one cusp
to one of the components of the complement of a closed geodesic in a
hyperbolic surface with cyclic fundamental group and no cusp.

For another interesting example, suggested by ~Bowditch, start with the
hyperbolic plane considered as the universal cover of the thrice
punctured sphere and the preimage of the three geodesics on the thrice
punctured sphere which run from one cusp to another, each with unit
transverse measure. Bend the plane. The closure $\xcl A$ of the
resulting convex surface meets $Q$ in the limit set of a Fuchsian
group which uniformizes a planar surface with three boundary circles,
together with sawteeth. That is, there is a spacelike plane $P$ such
that $\xcl A\cap Q\cap P$ is a Cantor set, and for each gap in the
Cantor set with endpoints $p,q$ a sawtooth consisting of a segment of
a line of the left ruling starting at $p$ and a segment of a line of
the right ruling starting at $q$, each going in the future direction
from $p$ or from $q$. The top surface of the convex hull determines a
left earthquake on the interior of the convex core of the planar
surface. The leaves of the shearing lamination spiral to the right
towards the boundary geodesics. After the earthquake they spiral to the
left towards the boundary geodesics. The earthquake with half as much
transverse measure turns the surface with geodesic boundary into a
thrice punctured sphere, an example of proposition ~\ref{p23}. Note
that the past boundary of the convex hull in anti de~Sitter space
contains parts of null planes. One expects the past boundary to
determine a right earthquake, defined on the projection to the
hyperbolic plane of the spacelike part of the boundary, and since the
spacelike part of the past boundary lies in a single plane, the
shearing lamination is the zero lamination. That is, the right
earthquake is the identity. Given any nowhere spacelike circle in $Q$
there is a corresponding left and right earthquake pair, each defined
on the interior of a convex subset of the hyperbolic plane with
geodesic boundary. Each conformally null segment of a line in the left
ruling in the boundary accounts for one geodesic on the frontier of
the range of the left earthquake and one geodesic on the frontier of
the domain of the right earthquake.

Let us close with some questions. Given a ``quasifuchsian'' group
acting on anti de~Sitter space, the volume of the convex hull and of
the domain of dependence are invariants. These volumes are functions
on the product of two copies of Teichm\"uller space. How do they
behave? Are they related, perhaps asymptotically, to such invariants
of a quasifuchsian group as the volume of the convex hull and the
Hausdorff dimension of the limit set? Suppose the hyperbolic structure
is specified on the two boundary components of the convex hull; does
there exist a unique manifold with the given pair of hyperbolic
structures? (In the case of quasifuchsian groups existence follows
from the Sullivan-Epstein-Marden theorem ~\cite{31,20} but I do not
know a proof of uniqueness.) Is a quasifuchsian group determined by
the hyperbolic structure on one boundary of the convex hull together
with the bending lamination on the other boundary? Is it determined by
the two bending laminations? Is it determined uniquely by the
conformal structure on the upper surface at infinity together with the
hyperbolic structure on the bottom surface of the convex hull?
Is it determined by the conformal structure of the upper surface at infinity together with the bending lamination on the bottom surface of the convex hull?
Laplace's equation can be solved by complexifying and separating
variables, so it looks like the wave equation. Analogously, one can
think of the left and right holonomy representations of a closed
spacelike surface in anti de~Sitter space as analogous to the two
conformal structures at infinity of a quasifuchsian group. So there
are analogous questions about the representation of the fundamental
group of a locally anti de~Sitter spacetime which is a neighbourhood of
a closed spacelike surface: Is it determined by the two measured
laminations, or by the hyperbolic structure on the future boundary of
the convex hull together with the measured lamination on the past
boundary of the convex hull, or by the left holonomy together with the
hyperbolic structure on one of the boundary components of the convex hull?

\section{Classification of spacetimes}\label{sec8}

Recall that a geometric structure on a manifold $M$ is (by definition)
complete provided that $M=X/\Gamma$ for some subgroup $\Gamma$ of $G$
acting properly discontinuously on $X$. We recall that $i)$ Margulis
~\cite{22,23}, cf. ~\cite{24,25} showed the existence of a complete
flat 3-dimensional Lorentz manifold with free fundamental group, $ii)$
complete closed Lorentz flat manifolds all have virtually polycyclic
fundamental group by a theorem of Goldman and Kamishima ~\cite{28} and
complete closed 3-dimensional affine flat manifolds all have
polycyclic fundamental groups by a theorem of Fried and Goldman
~\cite{26} and are classified in ~\cite{26}; the case of closed flat
Lorentz 3-manifolds having been done already in ~\cite{27}, and $iii)$
Carri\`ere ~\cite{5} shows that a compact affine manifold with linear
holonomy preserving a Lorentzian structure or more generally of
``discompacity one'' is complete. For more information see ~\cite{28}
and ~\cite{46}.

\begin{propo}\label{p24}
The linear holonomy of a complete flat spacetime is either solvable or else
a discrete subgroup of $\xO(2,1)$.
\end{propo}

\begin{proof}

Suppose $M$ is a complete oriented orthochronous spacetime. The
Lorentzian orthonormal oriented, time-oriented frame bundle of $M$ can
be identified with the quotient of the Lorentzian orthonormal oriented
time-oriented frame bundle of $\xr^{2+1}$ by the holonomy of $M$, and
the Lorentzian orthonormal oriented time-oriented frame bundle of
$\xr^{2+1}$ can be identified with $\iso(2,1)$. So $\pi_1M$ is a
discrete subgroup of $\iso(2,1)$. If $\pi_1M$ is not solvable and the
linear holonomy is not discrete, the linear holonomy must be a dense
subgroup of $\xso(2,1)_0$. Furthermore, by J{\o}rgensen's inequality, if
every $2$-generator subgroup of a finitely generated subgroup $H$ of
$\xso(2,1)_0$ is discrete, $H$ is discrete.

\begin{lemma}\label{l8}
Suppose $A,B$ generate a linear group which is not virtually solvable
and suppose that both $A,B$ have infinite order. Then for all
sufficiently large $N$, $A^N\text{ and }B^N$ freely generate a free group.
\end{lemma}

\begin{proof}
See ~\cite{30} and ~\cite{33}.
\end{proof}

We may assume the linear holonomy is torsion-free, passing if need be
to a finite index subgroup by Selberg's lemma. If the linear holonomy
is not discrete, there must be a subgroup $\langle A,B\rangle$
generated by two non-commuting elliptic elements. By Kronecker's
theorem on rational approximation, there is a sequence $N_i\to\infty$
such that $A^{N_i},B^{N_i}\to I$. By the previous lemma,  $\langle
A^{N_i},B^{N_i}\rangle$ is a free nonabelian group for all
sufficiently large $i$. So the holonomy contains a free group  
$\langle S_0,T_0\rangle$ such that
$S_0\mathbf{x}=A\mathbf{x}+\mathbf{c}$ and $T_0\mathbf{x}=B\mathbf{x}+\mathbf{d}$.
Inductively define
$S_{n+1}=[S_n,T_n],T_{n+1}=[T_n,S_n],A_{n+1}=[A_n,B_n],B_{n+1}=[B_n,A_n],S_n\mathbf{x}=A_n\mathbf{x}+\mathbf{c}_n,T_n\mathbf{x}=B_n\mathbf{x}+\mathbf{d}_n$.
Choose a norm $\|\quad\|$ on $\xr^{2+1}$. There is a constant $C$ such
that $\|I-[X,Y]\|+\|I-[Y,X]\|<C(\|X-I\|+\|Y-I\|)^2$ if
$\|X-I\|+\|Y-I\|$ is sufficiently small. We have
$[S_n,T_n]\mathbf{x}=A_{n+1}\mathbf{x}+A_{n}^{-1}C_{n}^{-1}(A-n\mathbf{d}_n-\mathbf{d}_n)+A_{n}^{-1}C_{n}^{-1}(\mathbf{c}_n-C_n\mathbf{c}_n)$
and similarly for $[T_n,S_n]$. So if $N_i$ is chosen sufficiently
large, $S_n,T_n$ converge to the identity, contradicting the fact that
the holonomy is discrete in $\iso(2,1)$.
\end{proof}

I believe proposition ~\ref{p24} is known.

\begin{propo}\label{p25}
There does not exist a complete flat spacetime $M$ with fundamental group
$\pi_1M\cong\pi_1S$ where $S$ is a closed surface of negative Euler characteristic.
\end{propo}

\begin{proof}
Suppose $M$ exists. Passing to a covering space we may assume $M$ is
oriented and orthochronous; then $S$ is an oriented surface. Since $M$
is complete, $\pi_1M$ injects into $\iso(2,1)$. Since $\pi_1M$ has no
normal abelian subgroup, the linear holonomy injects $\pi_1M$ into
$\xso(2,1)_0$. By proposition ~\ref{p24}, the image is discrete. So
$M$ has the holonomy of some standard spacetime.

$\widetilde{M}=\xr^{2+1}$ contains the universal covers of $M_1$ and
$M_2$, where $M_1,M_2$ respectively are the future, respectively the
past, of a future, respectively past directed closed convex spacelike
surface. $M_1$ and $M_2$ are disjoint by proposition ~\ref{p10}. Let
$M_0$ be the quotient of $\xr^{2+1}-\xint(M_1)\cup\xint(M_2)$ by
$\pi_1M$. Since the inclusion of either boundary component of $M_0$
into $M_0$ is a homotopy equivalence, $M_0$ is compact. By the
Thurston-Lok theorem, there is a family of flat Lorentz structures on
$M_0$ realizing all small deformations of the holonomy, and in this
family a neighbourhood of the original structure consists of spacetimes
with strictly convex spacelike boundaries. By the end of the proof of
proposition ~\ref{p3}, all these structures can be extended to
complete flat Lorentz structures on $M$. But even holding the linear
holonomy fixed, there are arbitrarily small deformations of the
holonomy which do not act fixed point free on $\xr^{2+1}$ and therefore cannot be the deck groups of a covering of $M$ by $\xr^{2+1}$. To see
this, observe that the measured geodesic laminations supported on
unions of simple closed curves are dense in the space of all measured
geodesic laminations, and by proposition ~\ref{p14}, or directly from
the construction of proposition ~\ref{p12}, the holonomy groups of such
structures cannot act freely on $\xr^{2+1}$. 
\end{proof}

Margulis ~\cite{22} conjectures that a complete flat Lorentz manifold with
finitely generated free fundamental group has discrete and purely
hyperbolic linear holonomy. By proposition ~\ref{p24}, the holonomy
must be discrete, but the possibility of parabolic elements is not
excluded. Propositions ~\ref{p24} and ~\ref{p25} show that a complete
flat Lorentz $2+1$ spacetime either has solvable fundamental group, in which
case the manifold is homotopy equivalent to a torus or to a Klein
bottle or else is a closed manifold, or else its fundamental group is
a discrete and non-cocompact subgroup of $\xO(2,1)$, in which case,
being torsion-free, it is a free group. It seems plausible that a
complete flat Lorentz manifold with free fundamental group is diffeomorphic
to the interior of a (possibly nonorientable) handlebody.

Propositions
~\ref{p6}, ~\ref{p7}, ~\ref{p8}, ~\ref{p15}, ~\ref{p26} motivate the
definition of the ``space of all classical solutions to Einstein's
equation in $2+1$ dimensions on a manifold of the form $S\times(0,1)$
containing $S$ as a closed spacelike surface'' as the family of all
domains of dependence of a fixed topological type. This family is
parametrized by the cotangent space of the Teichm\"uller space of $S$
and (cf. the remarks following proposition ~\ref{p13}) the cotangent
space has a natural symplectic structure. (In Witten's paper this
symplectic structure is obtained by reduction from an infinite
dimensional symplectic manifold.) There is a standard prescription for
``quantizing'' a cotangent bundle. Given a function on the cotangent
space which is either constant or linear on fibers, one can associate
an operator on the Hilbert space of square integrable half-densities on
Teichm\"uller space, namely a multiplication operator or a first order
differential operator. (A half density is a section of the square root
of the bundle of $n$-forms on a given manifold. For a half density
$g=f\sqrt{w},\|g\|^2=\int|f^2|w$ is invariantly defined on a manifold
with no given metric.) Equivalently one can use a Hilbert space of
functions on Teichm\"uller space, relative to the Weil-Petersson
volume element. In particular the trace of the linear holonomy of an
element of the fundamental group becomes an observable. In general the
holonomy of an element is conjugate to the composition of a hyperbolic
linear isometry of $\xr^{2+1}$ and a translation in the fixed
direction by a spacelike vector $\mathbf{v}$. The signed length
$l(\mathbf{v})$ of $\mathbf{v}$ is well-defined, using the convention
that the future pointing expanding eigenvector, the future pointing
contracting eigenvector, and an orthogonal unit spacelike vector
$\mathbf{w}$ form a positively oriented triad, and
$\mathbf{v}=l(\mathbf{v})\cdot\mathbf{w}$. The functions
$l(\mathbf{v})$ are linear on fibers, so they can be identified with
vector fields on Teichm\"uller space. I believe they are the
Hamiltonian vector fields associated to the traces of the linear
holonomy of group elements, but leave this question to the reader. The
Hilbert space together with these operators is called the solution of
quantum gravity in $2+1$ dimensions. In addition, Witten defines
remarkable new observables which do not come from the classical
observables, i.e., the invariants of the holonomy. (It may seem
strange to call a definition the solution of a problem, but until
Witten's work it was not clear how to make a reasonable definition.)

The classical uncertainty principle states that a wave function $\psi$
for which the variance of position
$\Delta_{x}^{2}=(x^2\psi,\psi)-(x\psi,\psi)^2$ is small has a large
uncertainty of momentum
$\Delta_{p}^{2}=((i\partial_x)^2\psi,\psi)-(i\partial_x\psi,\psi)^2$.
Because the bottom of the spectrum of the Laplacian on Euclidean space
is zero, there are wave functions $\psi$ for which the uncertainty of
momentum is arbitrarily small. The Weil-Petersson metric determines a
second order operator on Teichm\"uller space, which measures the
uncertainty of the translational part of the holonomy. Because the
mapping class group is non-amenable, we expect that the spectrum of the
Weil-Petersson Laplacian on Teichm\"uller space is bounded away from
zero. A covering of a compact Riemannian manifold with non-amenable deck
group has spectrum bounded away from zero, by Brooks's theorem
~\cite{45}; the moduli space is noncompact and indeed non-complete, so
the theorem does not apply, but still the result is likely to
hold. This means that for a wave function in the Hilbert space, there
is a lower bound on the uncertainty of the translational part of the
holonomy, and unlike the case of Euclidean space the bound is uniform:
the uncertainty cannot be made small by making the wave function very
widely spread over Teichm\"uller space. We leave these questions to
the reader.

In proposition ~\ref{p26} we give an extension of Carri\`ere's theorem
to flat Lorentz spacetimes with spacelike boundary, following his
proof as closely as possible for the reader's convenience. (Our proof
is terse, so the reader may wish to read Carri\`ere's proof first.) In
proposition ~\ref{p27} we will consider anti de~Sitter manifolds with
spacelike boundary.

\begin{definition}\label{d5}
Given a flat Lorentzian manifold (or more generally a Lorentzian
manifold of constant curvature) $N$ without boundary, let $E_x\subset
T_xN$ denote the set of vectors $v$ such that the geodesic flow $F_tv$
of $v$ is defined for all $t\in[0,1]$. $E_x$ is the domain of the exponential.
\end{definition}

\begin{propo}\label{p26}
Suppose $M$ is a compact flat $2+1$-dimensional orthochronous
Lorentzian manifold with spacelike boundary. Then if $\partial M$ is
empty $M$ is complete. Otherwise $M$ is a product
$\partial_0M\times[0,1]$ with each slice $\partial_0M\times\{t\}$ spacelike.
\end{propo}

\begin{proof}

First we enlarge $M$. If a component $\partial_{1i}M$ of the future
boundary is future directed there is a future complete manifold $F_iM$
with past boundary $\partial_{1i}M$. Adjoin each manifold $F_iM$ to
$M$ along $\partial_{1i}M$. If a component $\partial_{1i}M$ is past
directed, there is a manifold $G_iM$ such that the frontier of the
domain of dependence of the universal cover of $\partial_{1i}M$ is
locally convex surface $B_{1i}$ with null supporting planes, and
$G_iM$ is the quotient of the region between $B_{1i}$ and the
universal cover of $\partial_{1i}M$. Adjoin each manifold $G_iM$ to
$M$ along $\partial_{1i}M$. Similarly enlarge $M$ along the past
boundary. If there are toral boundary components adjoin their entire
futures or pasts (according as they are on the future or past
boundary.) Now we have an open manifold $M'$.

Let $(\widehat{M'},\;p:\widehat{M'}\to M')$ be the holonomy cover of
$M'$. $\widehat{M'}$ contains the holonomy cover $\widehat{M}$ of
$M$ as a closed submanifold with boundary. Given $x\in \widehat{M'}$
we will show that $E_x$, the domain of the exponential at $x$, is
convex. First we consider $F_x$, the intersection of $E_x$ with the
future pointing timelike and null vectors. Suppose $y,z\in F_x$ and
there is some $w\not\in E_x$ on the line segment $[y,z]$. Then there
exists $s_0$ maximal such that $[s_0y,s_0z]\subset F_x$ for $0\leq s\leq
s_0$. Consider the convex hull of the complement of $F_x$ in the cone
bounded by the rays through $y,z$. Because a dense subset of the
supporting planes of a convex set meet the convex set in an extreme
point, it is possible to replace $y,z$ by nearby points so that all of
the closed triangle $\triangle 0yz$ except one point $w$ on $[y,z]$
lies in $F_x$.

Now consider the geodesic ray $r=\xexp[0,w)$; it is the projection to
  $\widehat{M'}$ of the geodesic flow $F_t w$ of $w$ for $0\leq t <
  1$. If $r$ leaves the submanifold $\widehat{M}$ of  $\widehat{M'}$
  through a boundary component of  $\widehat{M}$ which is the past
  (respectively future) boundary of a future (respectively past)
  complete component of $\widehat{M'}-\widehat{M}$ then $r$ extends to
  a complete ray contradicting the fact that $\xexp[0,w)$ is a maximal
  ray. If $r$ leaves the submanifold $\widehat{M}$ and enters one of
  the preimages of a manifold $G_iM$, then the preimage in $E_x$ of the
  strictly convex surface $B_{1i}$ contains $w$ and is strictly convex
  toward $x$ (that is, an arc of the preimage of $B_{1i}$ together
  with rays from its endpoints to $0$ enclose a convex bounded set in
  the tangent space $T_x\widehat{M'}$.) This is impossible because all
  of $[y,z]$ except $w$ is in $E_x$. So the ray $r$ remains in
  $\widehat{M}$. Consider the projection $p(r)$ in $M$. Since $p(r)$
  is timelike, incomplete and maximal and $\partial M$ is compact and
  spacelike there is a neighbourhood $U$ of $\partial M$ such that at
  all sufficiently large times $r$, $p(r)$ is not in $U$. $M$ is
  compact so $p(r)$ is recurrent: There is an ellipsoid $\epsilon_0$ in the interior of $M$ such that $p(r)$ infinitely often enters the ellipsoid $p(\frac{\epsilon_0}{2})$
  where $\frac{\epsilon_0}{2}$ is the ellipsoid  obtained from $\epsilon_0$ by
  the affine map which fixes the center of $\epsilon_0$ and conjugates
  a translation by $v$ to a translation by $v/2$ for any
  $v\in\xr^{2+1}$. We suppose that $\epsilon_0$ is small enough that 
$p(\epsilon_0)$ is embedded. Let $\gamma_i\in\pi_1M$ be the
  elements defined by the successive entries of $p(r)$ into
  $\frac{1}{2}p(\epsilon_0)$. Now consider the preimages in the
  triangle $\triangle0yz$ by the exponential map of the ellipsoids
  $\gamma_i\cdot\epsilon_0$. After passing to a subsequence of
  $\{\gamma_i\}$, the preimages of the ellipsoids must converge in the
  Hausdorff topology (cf. ~\cite{5}) to a line segment or a line
  through $w$, and then the exponential of a line segment from $0$ to
  a point near $w$ on the segment $[y,z]$ must pass through
  $p(\epsilon_0)$ infinitely many times, which is impossible
  for a compact segment. So the domain $F_x$ is convex and similarly
  for the intersection $P_x$ of $E_x$ with the past pointing timelike
  and null vectors. The convergence of the subsequence of ellipsoids
  to a line or line segment is the crucial idea in Carri\`ere's proof.

Now suppose $E_x$ is not convex. As before we obtain $y,z$ in $E_x$
and a single point $w$ on $[y,z]$ which is not in $E_x$. If $y,z$ are
relatively timelike, we obtain a contradiction: The ray $p([0,w))$
cannot reach a preimage of a boundary component $B_{1i}$ because then
$[y,z]$ would map to a ray crossing $B_{1i}$. Now assume $y$ and $z$
are relatively spacelike. Suppose there are points $w_n$ arbitrarily
near $w$ on $[0,w)$ such that $p(\xexp(w_n))$ lies in $M$. Then
  (because each component of $\partial M$ is locally convex or locally
  concave) there is a slightly smaller compact manifold
  $M_0\subset\xint M$ such that $p(\xexp[0,w))$ is recurrent in
    $M_0$. Then there are ellipsoids $p(\gamma_n\cdot\frac{\epsilon_0}{2})$
    which cross the ray $p(\xexp[0,w))$ and (after taking a
      subsequence) the ellipsoids converge. The limit is necessarily
      the union of all the null lines through an interval of $[y,z]$
      containing $w$. (Otherwise some segment in the open triangle
      $\triangle0yz$ would cross infinitely many ellipsoids.) We
      obtain a contradiction because a compact segment $[0,w']$ where
      $w'$ is near $w$ on $[y,z]$ crosses infinitely many
      ellipsoids. So the ray $p(\xexp[0,w])$ eventually leaves
      $M\subset M'$. If it left through a concave boundary component --
      that is, a future, respectively a past boundary component of $M$
      which is the past, respectively future boundary of a future,
      respectively past complete end of $M$ then it eventually
      reenters $M$. But if it leaves through a convex boundary
      component we obtain a contradiction because the universal cover
      of the domain of dependence of a spacelike surface is convex, by
      proposition~\ref{p11}. So $E_x$ is convex.

By a result of Koszul (proposition 1.3.2 of ~\cite{5}),
$\widehat{M'}$ is $E_x$, regarded as an open subset of $\xr^{2+1}$. In
the case that $M$ is closed, Carri\`ere concludes that the frontier of
$\widehat{M}$ consists of $0,1\text{ or }2$ hyperplanes and deduces
that in fact $\widehat{M}=\xr^{2+1}$. Now suppose $\widehat{M}$ has at
least two future boundary components
$\partial_{1i}\widehat{M},\partial_{1j}\widehat{M}$. At most one of
these covers the past boundary in $M'$ of a future complete
submanifold of $M'$. Otherwise $\widehat{M'}$ would not be embedded in
$\xr^{2+1}$. Also at most one of
$\partial_{1i}\widehat{M},\partial_{1j}\widehat{M}$ is
convex. Otherwise a timelike linear functional would have two local
maxima on the closure of $\widehat{M'}$ in $\xr^{2+1}$. If one, say
$\partial_{1i}\widehat{M}$, is convex and another covers the past
boundary of a future complete manifold then it follows that the convex
manifold  $\widehat{M'}$ is all of $\xr^{2+1}$ contradicting the
assumption that $\partial_{1i}\widehat{M}$ is convex.

We conclude that if $M$ has a boundary component
$\partial_{i}M$ then $\pi_1\partial_{i}M$ has
index one in $\pi_1M$. If follows that (up to time reversal) either
$\widehat{M}$ has one future boundary component which is convex and
one which is concave, or else $M$ has two boundary components each of
which is a torus and $\widehat{M'}=\xr^{2+1}$. In the first case, for
any timelike direction $\mathbf{v}$ there is a family of planes
$P_t(\mathbf{v})=\{\mathbf{x}\cdot\mathbf{v}=ta(\mathbf{v})+(1-t)b(\mathbf{v})\}$
such that $P_0\mathbf{v}$ and $P_1\mathbf{v}$ are support planes of
the future and past boundaries respectively. For each $t$, the
envelope of all the planes $P_t$ is a spacelike surface, and this
surface defines an explicit product structure on $M$ with spacelike
slices. If the past or future boundary of $M$ is a torus whose
fundamental group has index one in $\pi_1M$ then the theorem holds
using proposition ~\ref{p8}.
\end{proof}

Let $X$ denote anti de~Sitter space and $\widetilde{X}$ the universal
cover of anti de~Sitter space. Suppose $M$ is a compact oriented
orthochronous anti de~Sitter manifold with spacelike boundary. We
consider $M$ as an $(\xiso\widetilde{X},\widetilde{X})$-manifold, so
the holonomy cover $\widehat{M}$ develops into $\widetilde{X}$.

\begin{propo}\label{p27}
a) Given a point $x\in\widehat{M}$, the exponential is a
diffeomorphism from the intersection of the domain $E_x$ of the
exponential with any spacelike plane in the tangent space
$T_x\widehat{M}$ to a convex set in a spacelike plane in
$\widehat{M}$. b) If $M$ is a compact oriented orthochronous anti
de~Sitter manifold with spacelike boundary, then either $M$ embeds in a
domain of dependence or else $M$ is a closed anti de~Sitter manifold
whose universal cover is a domain in $X$ with locally convex frontier. 
\end{propo}

\begin{proof}
The boundary components of $M$ must have negative Euler characteristic
by proposition ~\ref{p21}. Also, using proposition ~\ref{p21} $M$ can
be thickened to a manifold with locally strictly convex boundary. Let
$P$ be a spacelike plane in $\widetilde{X}$, and let $U$ be a
connected component of $\xdev^{-1}P$. Suppose $x\in U$. Let $T_xP$ be
the plane in $T_x\widehat{M}$ tangent to $U$. $P$ is a hyperbolic
plane, so the exponential map is injective on $T_xP$. The proof of
proposition ~\ref{p26} extends to show that $E_x$ meets $T_xP$ in a
convex set. Since $\widehat{M}$ has locally convex boundary, the
subset $\xexp(E_x\cap T_xP)\cap U$ is convex. Choosing a subgroup
$\xr=\widetilde{\xso(2)}\subset G_L$ defines a time function
$t:\widetilde{X}\to \xr$ such that the action of
$\xr$ permutes the planes $P_s=t^{-1}(s),s\in\xr$. So $\widehat{M}$ is
foliated by the preimages by $\xdev$ of the planes $P_s$. Now we will
show that the preimage of a spacelike plane $P$ is connected.

Suppose $\xdev^{-1}P_0$ is not connected, but contains disjoint
components $U,V$. There is some maximal interval $[a,b]$ containing
$0$ such that $\xdev^{-1}t^{-1}[a,b]$ contains two distinct components
$U',V'$ containing $U,V$. One but not both of $a,b$ may be infinite;
suppose $b$ is finite. Then $\xdev(U')\cap P_b$ and $\xdev(V')\cap P_b$ are
disjoint open convex sets, perhaps with disjoint closures, but for
$c>b$ sufficiently near $b$ there are points $p\in
U'\cap\xdev^{-1}P_b$, $q\in V'\cap\xdev^{-1}P_b$ and $x$ such that $\xdev(x)\in
P_c$ and such that $\widehat{M}$ contains a triangle $T$, i.e., a set
such that $\xdev(T)$ is the closed triangle
$\triangle\xdev(p)\xdev(q)\xdev(x)$ minus a closed interval in the
segment $\xdev(p)\xdev(q)$. The set $\triangle\xdev(p)\xdev(q)\xdev(x)$ lies
in a timelike plane. As in proposition ~\ref{p26} we obtain that
$\widehat{M}$ contains a triangle $\triangle p'q'x$ with one point
missing from the closed interval $p'q'$. By the argument of
proposition ~\ref{p26}, this is a contradiction.

So $\widehat{M}$ is foliated by open convex spacelike totally geodesic
surfaces and is embedded in $\widetilde{X}$. It follows that if
$\partial M$ is nonempty then $\widehat{M}$ has only one past
component and only one future component. Therefore the inclusion of
either boundary component of $M$ induces an isomorphism of fundamental
groups. Furthermore, $\widehat{M}$ has locally convex frontier in
$\widetilde{X}$ and at every point of the frontier there is a null
supporting plane. In this case either $\widehat{M}$ is contained in a
domain of dependence or else $\widehat{M}$ contains the entire region
bounded by the past boundary component $\partial_0\widehat{M}$ of
$\widehat{M}$ and a translate $T\cdot\partial_0\widehat{M}$ where $T$
generates the deck group of the covering of $X$ by $\widetilde{X}$ and
moves points to points in their futures. But then there would be a
domain of dependence $H(\partial_0\widehat{M})$ such that
$\pi_1\partial_0\widehat{M}$ acted properly discontinuously on a
region properly containing the entire future component of the frontier
of the domain of dependence $H(\partial_0\widehat{M})$ and this is
impossible. Indeed, the action on the future component of the frontier
of the domain of dependence is topologically equivalent to the action
of $\pi_1M'$ on the frontier of the universal cover of $M'$, for some
flat Lorentz manifold $M'$ which is a domain of dependence. So we may
assume $M$ is closed. We consider $M$ as a
$(\xiso\widetilde{X},\widetilde{X})$-manifold, so the development maps
the holonomy cover $\widehat{M}$ to $\widetilde{X}$. We have shown
that the development embeds $\widehat{M}$ as a simply connected open
subset of $\widetilde{X}$, since there is a submersion to $\xr$ whose
fibers are connected and simply connected. So we can regard $\pi_1M$
as a subgroup of $\xiso\widetilde{X}$. First we consider the case that
$\pi_1M$ meets the center $\xz$ of $\xiso\widetilde{X}$ in the
subgroup $n\xz,n>0$. In this case we can consider $M$ as an $(\xiso
X_n,X_n)$-manifold where $X_n$ is the $n$-fold cyclic cover of
$\xpsl_2\xr$; considered as such the holonomy cover (which we still
denoted by $\widehat{M}$) of $M$ develops to $X_n$ and this map is an
embedding such that the inclusion of the image is a homotopy
equivalence with $X_n$. We will identify $\widehat{M}$ with its
image. Fix $p\in\widehat{M}$. There is a smooth timelike loop
$c:S^1=\xr/\xz\to\widehat{M}$ which represents the generator of
$\pi_1X_n$. For each $t$ there is a timelike geodesic arc $d_t$,
parametrized by arc length, in $X_n$ joining $0$ to $c(t)$. Choose $c$
so that for $t\ne0\in S^1$, $c(t)$ is not one of the $n$ translates of
$p=c(0)$ by the center $\xz/n$ of $\xiso X_n$. Then $d_t$ depends
continuously on $t$. The set of $t$ for which $d_t$ lies in
$\widehat{M}$ is open. It is also closed. Indeed, suppose that for
$0\leq t\leq t_0$, $d_t$ lies in $\widehat{M}$. Let $L$ be the length
of $d_{t_0}$. The set of parameter values $s$ such that $d_{t_0}(s)$
is in $\widehat{M}$ is either all of $[0,L]$ or else there is some
$l<L$ such that $d_{t_0}([0,L))$ lies in $\widehat{M}$ but $d_{t_0}(L)$
  does not. Then after a small change in $p=c(0)$ and $d_{t_0}(L)$ we
  obtain, as in proposition ~\ref{p26}, a triangle $\triangle
  c(0)d_{t_0}(L)x$ in $X_n$ lying in a timelike plane such that all of
  the triangle except one point on the side $c(0)d_{t_0}(L)$ lies in
  $\widehat{M}$. By  Carri\`ere's recurrence argument (cf. proposition
  ~\ref{p26}) this is impossible. So $\widehat{M}$ contains a closed
  timelike curve through $p$, say with unit tangent vector
  $\mathbf{v}$ at $p$. The set of unit timelike vectors which
  are tangent to closed curves in $\widehat{M}$ is open, because all
  timelike geodesics in $X_n$ are closed and $\widehat{M}$ is open in
  $X_n$. It is also closed, by another application of Carri\`ere's
  argument. So $\widehat{M}$ contains the entire hyperbolic plane
  $\xh^2(p)$ which is one component of the preimage in $X_n$ of the
  dual plane in $\xpsl_2\xr=X_1$ to the image of the point $p$ in
  $X_1$. Given any point of $\xh^2(p)$, one (namely the normal) and
  therefore all of the unit timelike vectors through that point are
  tangent to a closed geodesic in $\widehat{M}$. We conclude that $M$
  is complete. $M$ has ``finite level'' in the terminology of ~\cite{53}.

Now suppose that $\pi_1M$ has trivial intersection with the center
$\xz$ of $\xiso\widetilde{X}$. $\widehat{M}$ is an open subset 
of $\widetilde{X}$. Suppose for the moment that $\widehat{M}$ is a 
proper open subset. The frontier of $\widehat{M}$ is locally convex
and 
defined by null supporting planes. Furthermore, $\widehat{M}$ is 
disjoint from its translates by the center of $\xiso\widetilde{X}$. 
Otherwise the frontier of $\widehat{M}$ would meet the interior of a
translate 
of $\widehat{M}$. But $\pi_1M$ acts properly discontinuously on the 
interior of $\widehat{M}$ and all of its translates, but not at any
point 
of the frontier, because any point of the frontier is the endpoint of
a ray 
which has recurrent projection to $M$ and therefore is the
accumulation 
point of an orbit. So we can consider $M$ as an $(\xiso X,X)$-manifold 
which is the quotient of a domain by a cocompact subgroup. 
The past and future boundaries of this domain are well-defined. 
By the arguments of  proposition ~\ref{p20} the frontier of either 
boundary component $\widehat{M}$ in $\xrp^3$ meets the quadric at
infinity 
in a nowhere timelike topological circle. In fact this is the entire 
intersection of the frontier with the quadric at infinity. 
For if the circles obtained from the past and future boundary
components 
were not equal, the quadric would contain some open set which lay
between them. 
Every point in the open set would be a point of the frontier. 
Otherwise there would be a geodesic ray leaving $\widehat{M}$ without 
passing through either the past or future boundary. 
So if the frontier of the past and future boundary components were not
equal, the frontier would meet the quadric in a closed set 
invariant under $\pi_1M$ and containing an open set. 
This closed set would necessarily equal all of the quadric at
infinity, 
which contradicts the hypothesis that $\widehat{M}$ is a proper 
open subset of $X$. 
So the frontier of $\widehat{M}$ in $\xrp^3$ meets the quadric at
infinity 
in a nowhere timelike topological circle, which we call the limit
circle. 
This contains at most countably many straight segments. 
If say the left holonomy had nontrivial kernel, any element in the 
kernel would fix all of the limit circle except possibly some open
segments. 
It would therefore fix set-wise three distinct null rays of the right 
ruling on the quadric and therefore also have trivial right holonomy. 
If follows that neither the left nor the right holonomies can have
nontrivial kernel. 
Suppose the limit circle contained segments of lines. 
Then one, say the left, holonomy leaves invariant a collection of open 
intervals in $\xrp^1$, from which it follows that the left holonomy
has 
free image. But since the left holonomy has trivial kernel, 
this implies that $M$ has free fundamental group. 
But $M$ is a closed aspherical manifold of dimension 3, so this is
impossible. 
So the limit circle contains no segments of lines. 
It therefore conjugates the left and right holonomies. 
Since $\xrp_{L}^{1}$ and $\xrp_{R}^{1}$ can be identified by the 
conjugating homeomorphism $\phi$ it follows from the discreteness of
$\pi_1M$ 
in $G_L\times G_R$, and the fact that $G_L$ can be identified with the 
space of distinct triples of points in $\xrp_{L}^{1}$ that $\pi_1M$ 
is discrete regarded as a subset of $G_L$.  But this implies that $\pi_1M$ is a Fuchsian group which contradicts the fact that $\pi_1M$ 
is the fundamental group of a closed aspherical 3-manifold. 
We conclude that $M$ is complete.

\end{proof}

Kulkarni and Raymond ~\cite{53} showed that a closed 3-manifold with a
complete anti de~Sitter structure is necessarily a Seifert
manifold. See also ~\cite{54} for more information on closed
3-manifolds with complete anti de~Sitter structure. ~\cite{54} shows
that given a complete Lorentzian manifold of the form
$\xpsl_2\xr/\Gamma$ where $\Gamma$ is a cocompact subgroup acting by
right translation, the deformations with left holonomy nontrivial but
lying in an abelian subgroup are still complete. By the Thurston-Lok
holonomy theorem, there are also deformations of the Lorentzian
structure corresponding to any small change in the holonomy; in
general the left holonomy will be irreducible after such a
change. Goldman conjectured in ~\cite{54} that these manifolds would
still be complete; this is established by proposition ~\ref{p27}. It
seems likely that by similar arguments one could prove there does not
exist a closed de~Sitter manifold and that a compact de~Sitter
manifold with spacelike boundary is a product in which all the slices
are spacelike, if one knew that a de~Sitter manifold which was a
tubular neighbourhood of a closed spacelike surface was associated to a
projective structure.

In the case of flat and anti de~Sitter spacetimes containing a closed
spacelike surface, we showed that a neighbourhood of the universal
cover of the surface embeds in the model space. The identity component
of the isometry group of the model space acts simply transitively on
the oriented orthochronous Lorentzian orthonormal frame bundle of the
model space. Since the holonomy acts properly discontinuously on the
neighbourhood of the surface, and therefore also on the restriction of
the frame bundle to this neighbourhood, the holonomy of the surface
must be discrete in the isometry group. This argument together with
proposition ~\ref{p24} is an alternative, in the case of a flat
spacetime, to the use of Goldman's theorem, and works in $n+1$
dimensions. Proposition ~\ref{p4} also generalizes. One can deduce that a
$3+1$-dimensional spacetime which is a small neighbourhood of a closed
spacelike hypersurface and has irreducible linear holonomy is a
deformation, by a 1-cocycle with values in $\xr^{3+1}$, of the
quotient of the positive (or negative) light cone by a group acting
cocompactly on the hyperbolic 3-space of unit timelike directions.

\end{document}